\documentclass[10pt,reqno,english]{amsart}
\usepackage{amsfonts,amsmath,latexsym,verbatim,amscd,mathrsfs,color,array}
\usepackage{amssymb}
\usepackage{pdfsync}
\usepackage{epstopdf}

\newcommand{\ass}{\quad\mbox{as}\quad}

\long\def\hide#1{}

\newcommand{\inn}{{\quad\hbox{in } }}

\newcommand{\ttt}{\tilde }

\newcommand{\LL}{{\tt L}  }
\newcommand{\TT}{{\mathcal T}  }

\newcommand{\nn}{ {\nabla}  }

\newcommand{\pp}{ {\partial} }

\newcommand{\vp}{\varphi}

\newcommand{\RR}{{{\mathbb R}}}

\newcommand{\GG}{ {\mathcal G}}

\newcommand{\C}{{\mathbb C}}
\newcommand{\R} {\mathbb R}
\newcommand{\Z} {\mathbb Z}
\newcommand{\F} {\mathbb F}
\newcommand{\SB} {\mathbb S}
\newcommand{\cuad}{{\sqcap\kern-.68em\sqcup}}

\newcommand{\dist}{{\rm dist}\, }
\newcommand{\foral}{\quad\mbox{for all}\quad}
\newcommand{\ve}{\varepsilon}

\newcommand{\be}{\begin{equation}}
	\newcommand{\ee}{\end{equation}}

\newcommand{\equ}[1]{(\ref{#1})}

\renewcommand{\div}{\mathop{\rm div}}
\newcommand{\curl}{\mathop{\rm curl}}

\newtheorem{lemma}{Lemma}[section]
\newtheorem{prop}{Proposition}[section]
\newtheorem{theorem}{Theorem}

\newtheorem{remark}{Remark}[section]
\newcommand{\bremark}{\begin{remark} \em}
	\newcommand{\eremark}{\end{remark} }

\numberwithin{equation}{section}

\begin{document}

	\title[Cluster of vortex helices  ]{Cluster of vortex helices in the incompressible $3$d Euler equations }

	\author[I. Guerra]{Ignacio Guerra}
	\address{\noindent I.G.: Universidad de Santiago de Chile (USACH) \\ 
		Facultad de Ciencia, Departamento de Matem\'atica y Ciencia de la Computaci\'on \\ 
		Las Sophoras 173, Estaci\'on Central, Santiago, Chile.}

	\email{ignacio.guerra@usach.cl}

	\author[M. Musso]{Monica Musso}
	\address{\noindent M.M.:  Department of Mathematical Sciences University of Bath,
		Bath BA2 7AY, United Kingdom.}
	\email{m.musso@bath.ac.uk}

	\maketitle
	
	\begin{abstract}
		In an inviscid and incompressible fluid in dimension $3$, we prove the existence of 
		several helical filaments, or vortex helices, collapsing into each others.
		
	\end{abstract}

	\section{Introduction}

	Ideal incompressible homogeneous fluid of density $\rho$ in three dimensional space in a time  interval $(0,T)$ 
	is governed by the Euler equations
	\be \label{euler0}
	\begin{aligned}
		{\bf u}_t  + ({\bf u}\cdot \nn ){\bf u} &= -{1\over \rho} \nn p  &&\inn \R^3 \times [0,T) , \\
		{\rm div}\, {\bf u} &= 0  && \inn \R^3\times [0,T)  ,
	\end{aligned}
	\ee
	supplemented by the initial velocity ${\bf u} (\cdot , 0) = u_0(x)$. Here $ {\bf u} \, : \,  \R^3 \times [0, T ) \to \R^3$ is the velocity field and $p \, : \, \R^3 \times [0, T ) \to \R$
	is the pressure, determined by the incompressibility condition. We will consider constant density
	$\rho = 1$.
	For a solution  ${\bf u }$ of  \equ{euler0},
	its vorticity is defined as
	$ \vec\omega = \curl {\bf u}  $. Then $\vec\omega$ solves the Euler system in {\em vorticity form}  \equ{euler0},
	\be \label{euler}
	\begin{aligned}
		\vec \omega_t  +
		({\bf u}\cdot \nn ){\vec \omega}
		&=( \vec \omega \cdot \nn ) {\bf u}  &&  \inn \R^3\times (0,T), \\
		\quad {\bf u}  = \curl \vec \psi,\ &
		-\Delta \vec \psi =  \vec \omega  && \inn \R^3\times (0,T),\\
		{\vec \omega}(\cdot,0)  &=  \curl u_0 && \inn \R^3 .
	\end{aligned}
	\ee
	Vortex filaments 
	are solutions of the Euler equations whose vorticity is concentrated in a small tube near an evolving imaginary smooth curve $\Gamma(t)$ embedded in entire $\R^3$  so that the associated velocity field vanishes as the distance to the curve goes to infinity. 
	Da Rios  \cite{darios}  in 1906, and  Levi-Civita \cite{levicivita1908} in 1908, found that, if  the vorticity
	concentrates smoothly and symmetrically in a small tube of size $\ve >0$ around a smooth curve for a certain interval of time,
	then it is possible to compute the instantaneous velocity of the curve to leading order.
	These computations suggest that the curve should evolve by its binormal flow, with a large velocity of order $|\log \ve |$. If $\Gamma(t)$ is parametrized as  $ \gamma(s,t)$ where $s$ designates its arclength parameter, then $\gamma(s,t)$ asymptotically obeys a law of the form
	$$
	\gamma_t =  2\, \kappa \, |\log\ve| \,   (\gamma_s\times \gamma_{ss})
	$$
	as $\ve \to 0$, or scaling $t= |\log\ve|^{-1}\tau $,
	\be\label{bin1}
	\gamma_\tau =  2\, \kappa \, (\gamma_s\times \gamma_{ss})= 2\, \kappa \, c \, {\bf b}_{\Gamma(\tau)}  .
	\ee
	Here $\kappa$ is the {\em circulation} of the velocity field on the boundary of sections to the filament, which is assumed to be a constant independent of $\ve$. Besides  ${\bf t}_{\Gamma(\tau)},\, {\bf n}_{\Gamma(\tau)},\, {\bf b}_{\Gamma(\tau)}$ are the usual tangent, normal and binormal unit vectors to $\Gamma (\tau )$, and $c$ its curvature. See \cite{ricca,majda-bertozzi} for a complete discussion on this topic.
	
	\medskip
	
	Jerrard and Seis \cite{jerrard-seis} rigorously proved  Da Rios' formal  computation,
	conditional upon knowing that the vorticity of a solution remains concentrated around some curve. In particular they consider a solution to \equ{euler} whose vorticity 
	$\vec \omega_\ve(x,t)$ satisfies, as $\ve \to 0$
	\be\label{conv}
	\vec \omega _\ve (\cdot,|\log\ve|^{-1}\tau)- \delta_{\Gamma(\tau)}  {\bf t}_{\Gamma(\tau)} \rightharpoonup  0,  \quad 0\le \tau \le T ,
	\ee
	where $\Gamma(\tau)$ is a sufficiently regular curve
	and $\delta_{\Gamma(\tau)}$ denotes a uniform Dirac measure on the curve. They proved that in these circumstances the curve
	$\Gamma(\tau)$ does indeed evolve by the law \equ{bin1}.  See \cite{jerrard-smets} and its references for results on the flow \equ{bin1}.

	\medskip
	The {\em vortex filament conjecture} (\cite{bcp}, \cite{jerrard-seis}) refers to the question of  existence of true solutions of 
	\equ{euler} that satisfy \equ{conv}  near a given curve $\Gamma(\tau)$
	that evolves by the binormal flow \equ{bin1}. This is an  open question, 
	except for very special cases. 
	
	\medskip
	A known solution of the binormal flow \equ{bin1} that does not change its form in time is a circle $\Gamma(\tau)$ with radius $R$
	translating with constant speed equal to ${2 \over R} $ along its axis.
	Solutions to \eqref{euler} whose vorticity is concentrated in a circular vortex filament are known as vortex rings, and the study of these objects dates back to Helmoltz and Kelvin' work. In 1894, Hill found an explicit axially symmetric solution of \equ{euler} supported in a sphere (Hill's vortex ring).
	Fraenkel's result \cite{fraenkel}
	provided the first construction of a vortex ring concentrated around  a torus with fixed radius and a small, nearly singular section $\ve>0$, travelling with constant speed $\sim |\log\ve |$, rigorously establishing the
	vortex filament conjecture for the case of traveling rings.  Vortex rings have been analyzed in larger generality in  \cite{fraenkel-berger,norbury,ambrosetti-struwe,devaleriola-vanschaftingen}.

	\medskip
	Another known solution of the binormal flow \equ{bin1} that does not change its form in time is the {\em rotating-translating helix}. It is a  circular helix of radius $R$ and pitch $h>0$, which rotates with constant speed 
	${ 2 \, \kappa \, h \over (h^2 + R^2)^{3\over 2}}$ and translates vertically in the direction of their axis of symmetry with constant speed ${2\, \kappa \, R^2 \over (h^2 + R^2 )^{3\over 2}}$.
	Solutions to the Euler equations \eqref{euler} whose vorticity is concentrated in a helical vortex filament are known as vortex helices, and the description of these objects started  with the works of  Joukowsky \cite{joukowsky}, Da Rios \cite{darios1916} and Levi-Civita \cite{levicivita1932}. In \cite{ddmw2}
	the authors provided the first construction of a vortex helix concentrated in an $\ve$-tubolar neighborhood of a rotating-translating helix evolving by binormal flow,  establishing the
	vortex filament conjecture for the case of rotating-translating helices. In \cite{ddmw2} the authors also find a solution to \eqref{euler} with several vortex helices, rotating-translating with comparable but different speeds.
	
	\hide{and the  affine transformation constituted by a rotation in the two first coordinates and a translation in the third one
		\begin{align*}
			S_{\theta, \sigma } ( x )   \, := \,   \left [    \begin{matrix}  P_\theta x' \\  x_3 + \sigma \end{matrix}    \right ]  =   Q_\theta x   + \left [   \begin{matrix}   0\\  \sigma  \end{matrix}    \right ], \quad x'=  \left [   \begin{matrix}   x_1\\ x_2 \end{matrix}    \right ].  
	\end{align*}}

	In this paper we are concerned with solutions to the Euler equations \eqref{euler} consisting of several vortex helices which are rotating-translating with almost the same speed. They are global-in-time solutions to \eqref{euler} and their vorticity has at main order the shape of several helical filaments which collapse into each others.  We call this phenomena a cluster of vortex helices. Let us be more precise.

	\medskip{}
	Let $N$ be a given integer. For any $i=1, \ldots , N$ consider points $(a_i,b_i)$ in $\R^2$, numbers $\sigma_i$ and $\beta_i$, and define the evolving curves $\Gamma_j $ parametrized by 
	\be \label{helix4}
	\gamma_i (s,\tau ) = \left( \begin{matrix}
		a_i \cos \big(
		\frac{s-\sigma_i \tau}{\sqrt{h^2 + R_i^2}} \big)-b_i \sin \big(
		\frac{s-\sigma_i \tau}{\sqrt{h^2 + R_i^2}} \big)
		\\
		a_i \sin \big(
		\frac{s-\sigma_i \tau}{\sqrt{h^2 + R_i^2}} \big) + b_i \cos \big(
		\frac{s-\sigma_i \tau}{\sqrt{h^2 + R_i^2}} \big)
		\\
		\frac{ h s + \beta_i \tau }{\sqrt{h^2+ R_i^2}}
	\end{matrix} \right) \in \R^3
	\ee
	where $h>0$ is a positive constant and  $$R_i = \sqrt{a_i^2 + b_i^2}.$$
	At any instant $\tau$, the curves $s \to \gamma_i (s,\tau)$ are  circular helices of radius 
	$R_i$ and common pitch $h$, parametrized by arc length. Their
	curvature  is $\frac{R_i}{R_i^2+h^2}$ and their torsion $\frac{h}{R_i^2 + h^2}$.  When time evolves, the curves rotate with constant speed $\frac{\sigma_i}{\sqrt{h^2 + R_i^2}}$ around the $z$-axis in $\R^3$ and at the same time translate vertically with constant speed $\frac{\beta_i}{\sqrt{h^2 + R_i^2}}$. At time $\tau =0$ each curve $\gamma_i$ passes through the point $(a_i, b_i, 0)$ in $\R^3$ (take $s=0$). A direct computation gives that $\gamma_i (s,\tau)$ evolves by the binormal flow \eqref{bin1} with circulation $\kappa_i$ provided the speeds $\sigma_i$ and $\beta_i$ are chosen to be 
	$$
	\sigma_i = \frac{2\, \kappa_i \, h}{R_i^2 + h^2}
	\quad 
	\beta_i = \frac{2\, \kappa_i\, R_i^2}{R_i^2 + h^2}.$$
	Each $\gamma_i(s,\tau)$  can be recovered from  $\gamma_i(s,0)$ by a rotation and a vertical translation
	\be\label{uff}
	\begin{aligned}
		\gamma_i (s,\tau)&=
		Q_{- {\sigma_i \over \sqrt{R_i^2 + h^2}}  \tau} \gamma(s,0)  + \left[ \begin{matrix}
			0 \\  {\beta_i \over \sqrt{R_i^2 + h^2}} \tau
		\end{matrix} \right], \quad \\Q_\theta &= \left( \begin{matrix} \cos \theta & -\sin \theta \\\sin \theta & \cos \theta\end{matrix} \right),
	\end{aligned}
	\ee
	or equivalently from $\gamma_i (0,0) = (a_i , b_i, 0)$ by mean of
	$$
	\gamma_i (s,\tau)=
	Q_{{ s-\sigma_i \tau \over \sqrt{R_i^2 + h^2}}  } (a_i , b_i, 0)   + \left[ \begin{matrix}
		0 \\  {s+\beta_i \tau\over \sqrt{R_i^2 + h^2}} 
	\end{matrix} \right].
	$$
	We shall identify the base point $(a_i, b_i, 0)$ of each helix simply with $(a_i,b_i)$. Helical filaments with comparable but different speeds as in \cite{ddmw2} have vorticity with
	\be\label{conv2}
	\vec \omega _\ve (\cdot,|\log\ve|^{-1}\tau) -\sum_{i=1}^N \delta_{\Gamma_i(\tau)}  {\bf t}_{\Gamma_i (\tau)}  \rightharpoonup 0 ,  \quad 0\le \tau \le T , \quad \ass \ve \to 0
	\ee
	so that $\gamma_i (0,0) = (a_i , b_i)$ do satisfy
	$$
	\dist \left( (a_i , b_i )\, , \, (a_j , b_j ) \right) > \delta, \quad \ass \ve \to 0, \quad \foral i\not= j 
	$$
	for some fixed $\delta >0$, independent of $\ve$.
	
	\medskip
	\noindent
	We are interested in colliding helical filaments. Let $r_0>0$ be a fixed number and assume that the points $(a_i, b_i)$ have the form, for all $i=1, \ldots , N$
	\be \label{helix50}
	(a_i,b_i)=(r_0,0)+ Q_i, \quad {\mbox {with}} \quad |Q_i| \to 0 \ass \ve \to 0 .
	\ee
	Since  $(a_i , b_i) \to (r_0,0)$ as $\ve \to 0$, for all $i$, the evolving helices $\gamma_i $ in \eqref{helix4} shrink into each others as $\ve \to 0$. The purpose of this paper is to establish the existence of a solutions to \eqref{euler} whose vorticity satisfies \eqref{conv2}, with colliding helical filaments $\Gamma_i$ in the sense \eqref{helix50}.
	
	\medskip
	\noindent
	We find that the points $Q_i$ needs to converge to $0$ at a precise  rate in terms of $\ve$. Let us be more precise. Assume 
	\be \label{helix5}
	(a_i,b_i)=(r_0 +s ,0)+ {{\bf P}_i \over |\log \ve|}, \quad \mbox{for}\quad i:=1\ldots N,
	\ee
	as $\ve \to 0$, for some constant $s$ and points ${\bf P}_i$ 
	satisfying 
	$$
	|s| < \delta {\log |\log \ve | \over |\ln \ve |} , \quad \delta < |{\bf P}_i | < \delta^{-1}
	$$
	for some $\delta >0$ small, and independent of $\ve$.
	The points ${\bf P}_i$ are at a uniform distance $d$ (independent of $\ve$) one from each other 
	\begin{equation}\label{defd}
		d= \min_{i\not=j} |{\bf P}_i - {\bf P}_j| >0.
	\end{equation}
	and the set $\{{\bf P}_1, \ldots , {\bf P}_N\}$ is  symmetric with respect to their first component, in the sense that
	\be\label{symmetry}
	{\bf P}= (p_1, p_2)  \in \{{\bf P}_1, \ldots , {\bf P}_N\}  \iff   (p_1, -p_2)  \in \{{\bf P}_1, \ldots , {\bf P}_N\}.
	\ee
	Writing
	\be\label{ff0}
	\quad {\bf P}_i= \left(\,  \, { {\bf P}_{i,1}  \over \sqrt{h^2 + r_0^2}} ,{\bf P}_{i,2}\right),
	\ee
	the points ${\bf P}_i$ satisfy at main order the balancing equations, for $i=1,\ldots,N$,
	\be \label{sysabintro0}
	\begin{aligned}
		&
		\sum_{j\not= i}  \kappa_j {({\bf P}_{i,1}-{\bf P}_{j,1})\over |{\bf P}_i - {\bf P}_j|^2} =\left(   \kappa_i {h r_0 \over 2\sqrt{(h^2+ r_0^2)^3} }   \, 
		-\alpha  {h r_0 \over 4\sqrt{h^2+ r_0^2} } \right)  \\
		& \sum_{j\not= i}   \kappa_j {({\bf P}_{i,2}-{\bf P}_{j,2})\over |{\bf P}_i - {\bf P}_j|^2} =0,
	\end{aligned}
	\ee   
	\medskip 
	where $\alpha$ is the constant defined by
	\be \label{alphadef}
	\alpha = \frac{2}{h^2+r_0^2} \, \frac{\sum_{i=1}^N \kappa_i^2}{\sum_{i=1}^N \kappa_i} .
	\ee
	We
	say that  $({\bf P}_1, \ldots , {\bf P}_N)$ of the form \eqref{ff0} is a {\em non-degenerate} solution to \eqref{sysabintro0}, if the linearization of the system \eqref{sysabintro0}  has only one element in its kernel, the one originating from the symmetry assumption \eqref{symmetry}.
	We will make this definition more precise in Section \S \ref{otto}.
	
	\medskip
	We prove the following result
	
	\begin{theorem}\label{teo2} 
		Let $h >0$, $r_0 >0$, $\kappa_1 , \dots \kappa_N$ be given numbers.
		Suppose there exists a {\em non-degenerate} solution $({\bf P}_1^0 , \ldots , {\bf P}_N^0)$ of the form \eqref{ff0} to system \eqref{sysabintro0}, satisfying \eqref{symmetry} and \eqref{defd}. 
		Let $\Gamma_j(\tau)$ be the helices parametrized by equation $\equ{helix4}$, for $j=1,\ldots,N$, with $(a_i , b_i)$ given by \eqref{helix5}. Then
		there exist $s^* \in \R$, points ${\bf Q}_1, \ldots , {\bf Q}_N$ and  a smooth solution $\vec \omega_\ve (x,t)$ to $\equ{euler}$, defined for $t\in (-\infty, \infty)$, such that
		$$
		(a_i,b_i)=(r_0 +s^* ,0)+ {{\bf P}_i \over |\log \ve|}, \quad {\bf P}_i = {\bf P}_i^0 + {\bf Q}_i , \quad |s^*| , |{\bf Q}_i | \lesssim {\log |\log \ve | \over |\log \ve |}$$
		and for all $\tau$,
		$$
		\vec\omega_\ve (x, \tau|\log\ve|^{-1}) -  \sum_{j=1}^N \kappa_j \delta_{\Gamma_j(\tau)}{\bf t}_{\Gamma_j(\tau)} \rightharpoonup  0 \ass \ve \to 0.
		$$
	\end{theorem}

	\medskip
	\noindent
	Our construction takes advantage of the invariance under helical symmetry of the Euler equations as discussed in \cite{dutrifoy, ettinger-titi, bnl,jln,ddmw2,velasco,bp,bd}. This invariance and the assumption that  the velocity field ${\bf u}$ in \eqref{euler0} is orthogonal to the helical symmetry lines imply that solutions to Problem \eqref{euler} can be found solving a transport equation in $2$ dimensions. 
	For a point 
	$$
	x = (x_1, x_2, x_3) \in \R^3, \quad x= (x', x_3), \quad x' \in \R^2
	$$
	consider the scalar transport equation for  $w(x',t)$ 
	\begin{equation}\label{PB-00}
		\begin{aligned}
			\left\{
			\begin{aligned}
				w_t + \nabla^\perp \psi \cdot \nabla w &=0 && {\mbox {in}} \quad \R^2 \times (0, T)
				\\
				- {\mbox {div}} (K \nabla \psi) &= w && {\mbox {in}} \quad \R^2 \times (0, T),
			\end{aligned}
			\right.
		\end{aligned}
	\end{equation}
	where $(a,b)^\perp = (b,-a)$ and
	$K (x_1,x_2) $ is the matrix
	$$
	K(x_1 , x_2 ) = \frac{1}{h^2+x_1^2+x_2^2}
	\left(
	\begin{matrix}
		h^2+x_2^2 & -x_1  x_2\\
		-x_1 x_2 & h^2+x_1^2
	\end{matrix}
	\right).
	$$
	Then there exists a vector field ${\bf u} = (u_1, u_2, u_3) $ with helical symmetry 
	$$
	{\bf u} ( Q_{\theta} x' , x_3+ h\theta) = \left( \begin{matrix} Q_\theta (u_1, u_2) \\ u_3 + h \theta \end{matrix} \right) \quad \forall \theta\in \R , \quad \forall x = (x',x_3) \in \R^3 ,
	$$
	such that
	\begin{align}
		\label{omegaHelical}
		\vec \omega(x,t)  =  \frac{1}{h} w( Q_{-\frac {x_3}h} x',t)  \, \left (\begin{matrix}   Q_{\frac \pi 2} x' \\ h    \end{matrix}  \right)  , \quad x = (x',x_3)
	\end{align}
	satisfies the Euler equation \eqref{euler}. Here $Q_\theta$ is the rotation matrix in the plane $(x_1, x_2)$ as defined in \eqref{uff}.
	The derivation of \eqref{PB-00} can be found for instance in \cite{dutrifoy, ettinger-titi, ddmw2}.

	\medskip
	Rotating solutions to problem \eqref{PB-00} with constant speed $\alpha$ have the form
	\begin{equation}
		\label{ansaztRot}
		w (x', \tau) = W \left( Q_{ \alpha \tau } x' \right), \quad \psi (x',\tau) = \Psi \left( Q_{ \alpha \tau } x' \right) .
	\end{equation}
	Let $\tilde x = P_{\alpha \tau} x'$.
	In terms of $(W, \Psi)$, the second equation in \eqref{PB-00} becomes
	$$
	-\div_{\tilde x} ( K(\tilde x) \nabla_{\tilde x} \Psi) = W,
	$$
	and the first equation gets the form
	\begin{align}
		\label{eq00}
		\nabla_{\tilde x} W \cdot \nabla_{\tilde x}^\perp \left( \Psi - \alpha |\log \ve | \frac{|\tilde x|^2 }{ 2} \right) = 0 .
	\end{align}
	See \cite{ddmw2} for details.
	We now observe that if $W(\tilde x) = F(\Psi(\tilde x) - \alpha  |\log \ve |
	\frac{ |\tilde x|^2}{2} )$, for some function $F$, then automatically \eqref{eq00} holds.
	We conclude that if $\Psi$ is a solution to
	\begin{equation}\label{PB10}
		- \div \cdot (K(\tilde x) \nabla_{\tilde x} \Psi ) = F (\Psi - \frac{\alpha}{2} |\log \ve | |\tilde x|^2 ) \quad {\mbox {in}} \quad \R^2,
	\end{equation}
	for some function $F$, and $W$ is given by
	$$
	W(\tilde x)= F\left( \Psi(\tilde x) - \alpha  |\log \ve | \frac{|\tilde x|^2}{2} \right)
	$$
	then $ (w, \psi )$ defined by \eqref{ansaztRot} is a solution for \eqref{PB-00}.
	
	\medskip 
	\noindent
	We now notice that a solution to \eqref{PB10} such that
	$$
	F (\Psi - \frac{\alpha}{2} |\log \ve | |\tilde x|^2 ) - \, \kappa_i \, \delta_{(a_i,b_i)} \rightharpoonup 0, \ass \ve \to 0
	$$
	gives a solution $\vec \omega (x,t)$ to \eqref{euler} of the form \eqref{omegaHelical} with the property that
	$$
	\vec \omega (x,t) - \, \kappa_i \, \delta_{\Gamma_i } {\bf t}_{\Gamma_i} \rightharpoonup 0, \ass \ve \to 0
	$$
	and $\Gamma_i $ defined as in \eqref{helix4} with
	$$
	{\sigma_i \over \kappa_i } \to \alpha \, h , \quad {\beta_i \over \kappa_i } \to {2r_0^2 \over r_0^2 + h^2}, \quad \ass \ve \to 0.$$

	\medskip
	\noindent
	The proof of Theorem \ref{teo2} is reduced to finding a non-linear function $F$ and a solution $\Psi$ to \eqref{PB10} such that
	\be \label{W}
	W (x):=  F (\Psi - \frac{\alpha}{2} |\log \ve | |\tilde x|^2 ) \sim  \sum_{i=1}^N \kappa_j \, \delta_{(a_i,b_i)}, \quad \ass \ve \to 0,
	\ee
	where $(a_i , b_i ) \to (r_0 , 0)$ for all $i$.
	We  build such a solution by means of elliptic singular perturbation techniques. For $N=1$ we recover the result in \cite{ddmw2}. Desingularization of point vortices for the Euler equations in dimension $2$ has been treated in \cite{ddmw,marchioro-pulvirenti,smetsvan}.

	\medskip
	Solutions concentrated near helices in the Euler equations and also other PDE settings have
	have been built in \cite{cao1,cao2,contreras-jerrard,chiron,ddmr,jerrard-smets-gp,wei-yang}.  Solutions to the Euler equations \eqref{euler} concentrated around several vortex rings which are collapsing one into each others are known in literature. The first result is due to Buffoni \cite{buffoni}, who constructed co-axial vortex rings
	(sets homeomorphic to solid tori) moving along their common axis at the
	same propagation speed. These rings are nested in the sense that the convex-hull of one ring contains the subsequent ring and at the same time the two rings do not intersect. A more recent result is contained in \cite{ALW}, where they also find the formal law for the dynamics of the centers of a family of clustering rings. It happens that this is  the same law \eqref{sysabintro0} that governs the helical clustering phenomena.   The same law of motion also governs the interaction of multiple vortex rings in the Gross-Pitaevskii equation \cite{AHLW, jerrard-smets-gp}. The law for the interaction of nearly parallel vortex filaments has been studied in \cite{klein-majda-damodaran}.

	Configurations of points $({\bf P}_1^0, \ldots , {\bf P}_N^0)$ that satisfy \eqref{sysabintro0} and the assumptions of Theorem  \ref{teo2} are known in the literature. For instance, letting 
	$N= n+m$ and $\kappa_i=1$ for $i=1,\ldots,m$ and  $\kappa_i=-1$ for $i=m+1,\ldots,n+m$, explicit non-degenerate solutions to \eqref{sysabintro0} when  
	$$
	(m,n)\in \SB:=\{(2,1),(3,2),(4,3),(5,4),(6,5)\}
	$$
	are described in \cite{AHLW}.
	In this case
	a direct computation gives
	$$
	\alpha=\frac{2(m+n) }{(h^2+r_0^2)(m-n)}
	$$
	from which we deduce that $m$ must be different from $n$. Other constructions of admissible configurations can be found in  \cite{AHLW,ALW}.
	
	%
	%
	

	\medskip
	As we already discussed, Theorem \ref{teo2} follows from proving the existence of a function $\Psi$ and a non-linearity $F$ to solve \eqref{PB10} and \eqref{W}. This is what the rest of the paper is devoted to.
	In Section \S \ref{sec3} we find a smooth stream function $\Psi$ solving approximately 
	$$
	- \div \cdot (K( x) \nabla\Psi ) \sim \sum_{i=1}^N \kappa_i \delta_{(a_i,b_i)},
	$$
	in coherence with the expectation \eqref{W}. In Section \S \ref{sec4} we choose the non-linearity $F$. It will be reminiscent  of 
	$
	f(s) = e^s
	$
	and the Liouville equation $\Delta u + e^u=0$ in $\R^2$ will be used as a limit problem to describe the profile of the helical filaments, near the centres of the vortex helices (see \eqref{defF}-\eqref{defF1}). We define a first approximate solution in Section \S \ref{sec4}, and estimate the error of approximation in Section \S \ref{sec5}. After the approximate solution is built, we proceed to find an actual solution close to the approximation. The actual solution is found using the inner-outer gluing method, which has been used in several other contexts. References for problem related to inviscid incompressible fluids are \cite{ ddmw} for the problem of point vortex  desingularization for the Euler equations in dimension $2$, \cite{ddmw3} for the leapfrogging of vortex rings, and also \cite{ddmw2}. Since the interaction among different helices is strong (as their relative distance is small), it is relevant for us to pose the inner-outer gluing method so that these interactions can be controlled. Section \S \ref{appe} contains two basic elliptic linear theories which are at the core of the resolution of the inner-outer scheme. Sections \S \ref{sette} and \ref{otto} are devoted to find an actual solution to the problem, where the choice of the centers of the helices to solve \eqref{sysabintro0} plays a central role.  
	
	\section{Finding the approximate stream function}\label{sec3}
	Let $h >0$, $r_0 >0$, $\kappa_1 , \dots \kappa_N$ be given numbers, and define $\alpha $ to be the constant defined in \eqref{alphadef}. 
	
	The rest of the paper is devoted to find a solution $\Psi $ to the semi-linear elliptic equation
	\begin{equation}\label{PB1}
		\nabla \cdot (K \nabla \Psi ) + F(\Psi -\frac{\alpha}{2} |\log \ve | |x|^2 )= 0 \quad {\mbox {in}} \quad \R^2.
	\end{equation}
	More precisely   we look for a non-linear function $F$ and a solution $\Psi$ to \eqref{PB1} with the property that, if we set 
	$$W (x) = F(\Psi -\frac{\alpha}{2}|\log \ve | |x|^2 )
	,$$
	then
	\be \label{vort}
	W (x) \sim 8 \pi \sum_{j=1}^N \kappa_j \, \delta_{P_j}, \quad \ass \ve \to 0
	\ee
	for some  $P_j \in \R^2$. The points $P_j$ are assumed to be close to each other and collapse to the same point, as $\ve \to 0$. We assume they have the form 
	\be \label{points}
	\begin{aligned}
		P_j&= (r_0 + s, 0) + \frac{\hat P_j}{|{\log \ve }|}, \quad P= (P_1, \ldots , P_N), \quad \hat P= (\hat P_1, \ldots , \hat P_N).
	\end{aligned}
	\ee
	We assume  the following bounds on $s$ and $\hat P$
	\begin{equation}\label{points1}
		\begin{aligned}
			\| \hat P \| \lesssim  1 , \quad |s| \lesssim {\log |\log \ve | \over |\log \ve |}
		\end{aligned}
	\end{equation}
	In other words, we look for  the stream function $\Psi$ to have the asymptotic behaviour  
	$$
	\Psi (x) \sim \sum_{j=1}^N \kappa_j \Psi_j (x) , \quad {\mbox {with}} \quad -\nabla \cdot (K \nabla \Psi_j ) \sim 8 \pi \delta_{P_j}, \quad \ass \ve \to 0.
	$$
	We expect each function $\Psi_j$ to be, locally around $P_j$, an approximate Green's function for the operator $\nabla \cdot (K \nabla \cdot )$ in $\R^2$. 
	
	\medskip
	This section is devoted to analyze the approximate Green's function for the operator $\nabla \cdot (K \nabla \cdot )$ in $\R^2$ and to construct an approximate stream function for the $N$-helical filaments.

	\subsection{Approximate Green's function for the operator $\nabla \cdot (K \nabla \cdot )$ in $\R^2$.} The purpose of this sub-section is to find an explicit regular function which locally around a point $P= (a,b) \in \R^2$  satisfies approximately
	\be \label{g1}
	-\nabla \cdot \left( K \nabla \Psi \right) = 8 \pi \delta_P.
	\ee
	To this purpose we need to understand the structure of the operator in divergence form 
	\be \label{defL}
	\begin{aligned}
		L&:=-\nabla \cdot ( K \, \nabla \, ), \quad {\mbox {where}} \quad \\
		K&= \frac{1}{h^2+x_1^2+x_2^2}
		\left(
		\begin{matrix}
			h^2+x_2^2 & -x_1  x_2\\
			-x_1 x_2 & h^2+x_1^2
		\end{matrix}
		\right)
	\end{aligned}
	\ee
	when evaluated around a given point $P$. We will show that, after an ad-hoc change of variable, the operator $L$ will look like the usual Laplace operator in $\R^2$ when considered in a neighborhood of $P$.
	
	\medskip
	The operator $L$ is explicitly given by
	\begin{equation}\label{exL1}
		\begin{aligned}
			L& = {h^2+ x_2^2 \over h^2+ r^2} \pp_{x_1 x_1} + {h^2+ x_1^2 \over h^2+ r^2} \pp_{x_2 x_2} - 2{ x_1 x_2 \over h^2+r^2} \pp_{x_1 x_2}\\
			&- {x_1 \over h^2+ r^2} \left({2 h^2\over h^2+r^2} +1\right) \pp_{x_1} - {x_2 \over h^2+ r^2} \left({2 h^2\over h^2+r^2} +1\right) \pp_{x_2} .
		\end{aligned}
	\end{equation}
	Indeed,  using the notation $K = \left( \begin{matrix} K_{11} & K_{12} \\ K_{12} & K_{22} \end{matrix} \right)$, we get
	\begin{align*}
		L& = \nabla \cdot (K \nabla \cdot) = K_{11} \, \pp_{x_1}^2 + K_{22} \, \pp_{x_2}^2 + 2 K_{12} \,  \pp_{x_1 x_2}^2\\
		&+ ( \pp_{x_1} K_{11}   + \pp_{x_2} K_{12} ) \,  \pp_{x_1} + ( \pp_{x_2} K_{22} + \pp_{x_1} K_{12}) \, \pp_{x_2} \\
		&= {h^2+ x_2^2 \over h^2+ r^2} \pp_{x_1 x_1} + {h^2+ x_1^2 \over h^2+ r^2} \pp_{x_2 x_2} - 2{ x_1 x_2 \over h^2+r^2} \pp_{x_1 x_2}\\
		&+ \left( \pp_{x_1} \left({h^2+x_2^2 \over h^2+r^2} \right) - \pp_{x_2} \left({x_1 x_2 \over h^2+r^2} \right) \right) \pp_{x_1} \\
		&+ \left( \pp_{x_2} \left({h^2+x_1^2 \over h^2+r^2} \right) - \pp_{x_1} \left({x_1 x_2 \over h^2+r^2} \right) \right) \pp_{x_2} ,
	\end{align*}
	where $r=|x|$.
	Formula \eqref{exL1} follows directly from the facts that
	$$
	\pp_{x_1} \left({h^2+x_2^2 \over h^2+r^2} \right) - \pp_{x_2} \left({x_1 x_2 \over h^2+r^2} \right)= - {x_1 \over h^2+ r^2} \left({2 h^2\over h^2+r^2} +1\right)
	$$
	and
	$$
	\pp_{x_2} \left({h^2+x_1^2 \over h^2+r^2} \right) - \pp_{x_1} \left({x_1 x_2 \over h^2+r^2} \right)= - {x_2 \over h^2+ r^2}\left({2 h^2 \over h^2+r^2} +1\right).
	$$
	
	\medskip	
	Let us introduce the change of variables
	\begin{align*}
		x_1-a&={a  h \over R \sqrt{h^2+ R^2}  } z_1 - {b \over R} z_2 \\
		x_2 - b &= {b h  \over R \sqrt{h^2+ R^2}  } z_1 + {a \over R} z_2,
	\end{align*}
	where $
	R= \sqrt{a^2 + b^2}.
	$ This is equivalent to say
	\begin{align*}
		z_1&= { \sqrt{h^2+ R^2}  \over h R  } \left[ a (x_1 - a) + b (x_2 - b)\right] \\
		z_2 &= {1 \over R} \left[ -b (x_1 - a) + a (x_2 - b)\right].
	\end{align*}
	We will also use the  matrix notation
	\be \label{matrixdef}
	x-P = A[P] z, \quad A[P] = \left( \begin{matrix} 
		{a h\over R \sqrt{h^2+ R^2} } & - {b \over R}  \\
		{b h \over R \sqrt{h^2+ R^2} } & {a \over R} 
	\end{matrix} \right).
	\ee
	When expressed in the $z$-variable, we recognize that the operator $L$ takes the form
	\be \label{Lz}
	L= \Delta_z + B,
	\ee
	where
	\be \label{B0}
	\begin{aligned}
		B&=\left( {h^2(R^2-r^2)+z_2^2(h^2+R^2) \over (h^2+ r^2)h^2}\right) \pp_{z_1 z_1} \\
		&+ {1 \over (h^2+ r^2)}\left(\left(z_1\frac{h}{\sqrt{h^2+ R^2}}+R\right)^2-r^2\right) \pp_{z_2 z_2} \\
		&-2{\sqrt{h^2+R^2} \over h (h^2+ r^2)}z_2\left( z_1
		{h \over \sqrt{h^2+ R^2}} +R \right) \pp_{z_1 z_2}\\
		&- \frac{z_1(h^2+R^2)+Rh\sqrt{h^2+ R^2}}{h^2(h^2+r^2)} \left({2 h^2\over h^2+ r^2} +1 \right)  \pp_{z_1 }\\
		&
		- \frac{z_2}{h^2+r^2} \left({2 h^2\over h^2+ r^2} +1 \right) \pp_{z_2 }.
	\end{aligned}
	\ee
	Here $r=r(z)$ is
	\begin{align*}
		r^2=|x|^2 = R^2 + 2 R {h \over \sqrt{h^2+ R^2}} z_1 + q_2 (z) \:\:
		\mbox{with}\:\:
		q_2 (z) = { h^2 \over h^2+R^2 } z_1^2  + z_2^2.
	\end{align*}
	Formulas \eqref{Lz} and \eqref{B0} are consequence of the following  straightforward computations
	\begin{align*}
		\pp_{x_1} &= {a \over R} {\sqrt{h^2+ R^2} \over h} \pp_{z_1} -{b \over R} \pp_{z_2}\\
		\pp_{x_2} &= {b \over R} {\sqrt{h^2+ R^2} \over h}\pp_{z_1} +{a \over R} \pp_{z_2}\\
		\pp_{x_1 \, x_1} &= {a^2 \over R^2} \, {(h^2+ R^2) \over h^2}\, \pp_{z_1 \, z_1} + {b^2 \over R^2} \pp_{z_2 z_2} - 2 {ab \over R^2} {\sqrt{h^2+ R^2} \over h}\pp_{z_1 \, z_2}\\
		\pp_{x_2 \, x_2} &= {b^2 \over R^2} \, {(h^2+ R^2) \over h^2}\, \pp_{z_1 \, z_1} + {a^2 \over R^2} \pp_{z_2 z_2} + 2 {ab \over R^2} {\sqrt{h^2+ R^2} \over h}\pp_{z_1 \, z_2}\\
		\pp_{x_1 \, x_2} &= {ab \over R^2} \, {(h^2+ R^2) \over h^2}\, \pp_{z_1 \, z_1} - {ab\over R^2} \pp_{z_2 z_2} + {a^2 - b^2 \over R^2} {\sqrt{h^2+ R^2} \over h}\pp_{z_1 \, z_2}.\\
	\end{align*}
	
	\medskip
	The operator $B$ in \eqref{B0} becomes a small perturbation of the Laplacian, when we restrict our attention to a small region around the point $P$, that in the $z$-variable can be described with  $|z| < \delta$, for
	a fixed $\delta$ small. Indeed, in this region  the operator $B$ has  the form
	$$
	\begin{aligned}
		B&=\left( -2{R h \over  (h^2+ R^2)^{3/2}} \, z_1 + O(|z|^2 )\right) \pp_{z_1 z_1} + O(|z|^2 ) \pp_{z_2 z_2} \\
		&-\left( 2{R \over h\sqrt{h^2+ R^2} } z_2 + O(|z|^2) \right) \pp_{z_1 z_2}
		\\
		&
		- \left( {R \over h\sqrt{h^2+ R^2} } \left({2h^2\over h^2+ R^2} +1 \right)  + O(|z|) \right) \pp_{z_1 } \\
		&- \left(\frac{z_2}{h^2+R^2} \left({2 h^2\over h^2+ R^2} +1 \right)+O(|z|^2)\right) \pp_{z_2 }.
	\end{aligned}
	$$
	Equation \eqref{g1} thus becomes
	\be \label{g2}
	-\left( \Delta + B \right) \psi = 8 \pi \delta_0, \quad \psi (z ) = \Psi \left( P + A[P] z \right).
	\ee
	We now choose a regularization of the Green's function of the Laplace operator $\Delta_z$ as a starting point for the construction of the approximate regularization to \eqref{g2}. The regularization we choose is a radial solution of the Liouville equation
	\be \label{liouville}
	\Delta u + e^u = 0 \quad \inn \R^2, \quad \int_{\R^2} e^u \, < \infty.
	\ee
	All solutions to \eqref{liouville} that are radially symmetric with respect to the origin are given by 
	$$
	\Gamma_{\mu \ve} (z) - 2 \log \ve \mu , \quad {\mbox {where}} \quad \Gamma_{\ve  \mu } (z) = \log {8 \over (\ve^2 \mu^2 + |z|^2 )^2}
	$$
	for any value of the constants $\ve$ and $\mu>0$. Indeed we have
	$$
	-\Delta \Gamma_{\ve  \mu }=\ve^2 \mu^2 e^{\Gamma_{\ve \mu}} = {1\over \ve^2 \mu^2} U \left( {z \over \ve \mu} \right)  \quad {\mbox {with}} \quad U(y) = \frac8{(1+ |y|^2)^2}.  
	$$
	Hence
	$$
	-\Delta \Gamma_{\ve  \mu } \rightharpoonup 8 \pi \delta_0 , \quad \ass \ve \mu \to 0.
	$$
	Besides, we compute
	\begin{align*}
		B [ \Gamma_{\ve  \mu }] &= - {2 Rh \over (h^2+R^2)^{3/2} } z_1 \pp_{z_1 z_1} \Gamma_{\ve  \mu }  -  {2 R \over h\sqrt{h^2+ R^2}}  z_2 \pp_{z_1 z_2} \Gamma_{\ve  \mu }  \\
		&- {R \over h\sqrt{h^2+ R^2}} \left(1+ {2h^2 \over h^2+ R^2} \right) \pp_{z_1} \Gamma_{\ve  \mu }+ E_1
	\end{align*}
	where $E_1$ is a smooth  function, uniformly bounded for $\ve \mu $ small, in a bounded region for $z$.
	
	\medskip
	We take advantage of the  explicit expression of $\Gamma_{\ve  \mu }$ to  find
	\begin{align*}
		\pp_{z_1} \Gamma_{\ve  \mu } (z) = -{4z_1 \over \ve^2 \mu^2+|z|^2} , &\quad  z_1 \pp_{z_1z_1} \Gamma_{\ve  \mu } (z) = -{4 z_1\over \ve^2 \mu^2 +|z|^2} +{8z_1^3 \over (\ve^2  \mu^2 +|z|^2)^2}\\
		z_2 \pp_{z_1z_2} \Gamma_{\ve  \mu } (z) &= {8z_2^2 z_1 \over (\ve^2 \mu^2 +|z|^2)^2}.
	\end{align*}
	Using that
	$$
	z_1 z_2^2  = { |z|^2 z_1 \over 4} - { {\mbox {Re} } (z^3 )\over 4},\quad z_1^3 = { 3|z|^2 z_1 \over 4} + { {\mbox {Re} } (z^3 )\over 4}
	$$
	we obtain
	\begin{align*}
		-&{2 R h \over (h^2+R^2)^{3/2} } z_1 \pp_{z_1 z_1} \Gamma_{\ve  \mu } -  {2 R \over h\sqrt{h^2+R^2}} z_2 \pp_{z_1 z_2} \Gamma_{\ve  \mu } \\
		&- {R \over h\sqrt{h^2+R^2}} \left(1+ {2 h^2\over h^2+ R^2} \right) \pp_{z_1} \Gamma_{\ve  \mu } \\
		&  = \left[{8h R \over (h^2+ R^2)^{3\over 2} } + {4R(3h^2+ R^2) \over h (h^2+ R^2)^{3\over 2} }\right]{z_1 \over \ve^2 \mu^2 + |z|^2}   \\
		&- {16 h R  \over (h^2+ R^2)^{3\over 2} }{z_1^3 \over (\ve^2 \mu^2 + |z|^2)^2}
		- {16 R  \over h(h^2+ R^2)^{1\over 2} }{z_1 z_2^2 \over (\ve^2 \mu^2 + |z|^2)^2} \\
		&= {4h R \over (h^2+ R^2)^{3\over 2} } {z_1 \over \ve^2 \mu^2 + |z|^2} + {4 R^3 \over h (h^2+ R^2)^{3\over 2}} {{\mbox {Re} } (z^3 ) \over (\ve^2 \mu^2 + |z|^2)^2} \\
		&+ {4 R (4h^2+ R^2) \over h (h^2+ R^2)^{3\over 2} }  { \ve^2 \mu^2 \, z_1 \over (\ve^2 \mu^2 + |z|^2 )^2}
	\end{align*}
	We can modify the function $\Gamma_{\ve \mu}$ to eliminate  part of the above error. We define
	\be \label{defc}
	\begin{aligned} 
		\psi_{1} (z) &= \Gamma_{\ve \mu} (z) \left( 1+ c_1 z_1  \right), \quad {\mbox {with}} \\
		c_1 & ={1\over 2}  {R h \over (h^2+R^2)^{3\over 2} }\\
	\end{aligned}
	\ee
	Since
	$$
	(\pp_{z_1 z_1} + \pp_{z_2 z_2} ) (c_1 z_1 \Gamma_{\ve \mu} ) = - 8 \, c_1 \, {z_1 \over \ve^2 \mu^2 + |z|^2} - 8\,  c_1 \, {\ve^2 \mu^2 z_1 \over (\ve^2 \mu^2 + |z|^2 )^2}, 
	$$
	with this choice of  $c_1$
	we get
	\begin{align*}
		L (\psi_1 ) (z)&= \Delta \Gamma_{\ve \mu} +  {4 R^3 \over h(h^2+ R^2)^{3\over 2}} {{\mbox {Re} } (z^3 ) \over (\ve^2 \mu^2 + |z|^2)^2} + {4 R (3h^2+ R^2) \over (h^2+ R^2)^{3\over 2} }  { \ve^2 \mu^2 \, z_1 \over (\ve^2 \mu^2 + |z|^2 )^2}\\
		& 
		+ c_1 B(z_1 \Gamma_{\ve \mu} ) 
		+ E_1
	\end{align*}
	where $E_1$ is an explicit function, which smooth in the variable $z $ and uniformly bounded, as $\ve \mu \to 0$.
	
	\medskip
	We now use the fact that 
	$z_2\pp_{z_2} \Gamma_{\ve  \mu } (z) = -\dfrac{4z_2^2}{ \ve^2 \mu^2+|z|^2} $, to write
	\begin{align*}
		c_1 B(z_1 \Gamma_{\ve \mu} ) &= c_1z_1B(\Gamma_{\ve \mu} ) - c_1{4 Rh \over (h^2+R^2)^{3/2} } z_1 \pp_{z_1} \Gamma_{\ve  \mu }  -  c_1{2R \over h\sqrt{h^2+ R^2}}  z_2 \pp_{z_2} \Gamma_{\ve  \mu }  \\
		&- c_1{R \over h\sqrt{h^2+ R^2}} \left(1+ {2h^2 \over h^2+ R^2} \right)  \Gamma_{\ve  \mu }
		+ \bar E_1 \\
		&= - {R^2 \over 2(h^2+ R^2)^2 } \left({2h^2\over h^2+ R^2} +1 \right)  \Gamma_{\ve \mu}
		+ E_2 
	\end{align*}
	where $E_2$ is another explicit function,  smooth in the variable $z $ and uniformly bounded, as $\ve \mu \to 0$.
	
	\medskip
	Combining these computations we obtain that the function $\psi_1$ introduced in \eqref{defc} satisfies
	\begin{align*}
		L (\psi_1 ) (z)&= \Delta \Gamma_{\ve \mu} +  {4 R^3 \over h(h^2+ R^2)^{3\over 2}} {{\mbox {Re} } (z^3 ) \over (\ve^2 \mu^2 + |z|^2)^2} + {4 R (3h^2+ R^2) \over h(h^2+ R^2)^{3\over 2} }  { \ve^2 \mu^2 \, z_1 \over (\ve^2 \mu^2 + |z|^2 )^2}\\
		&- {R^2 \over 2(h^2+ R^2)^2 } \left({2h^2\over h^2+ R^2} +1 \right)  
		\Gamma_{\ve \mu} + E_1 + E_2 ,
	\end{align*}
	where $E_1$ and  $E_2$ are explicit functions, smooth in the variable $z $ and uniformly bounded, as $\ve \mu \to 0$.

	\medskip
	Our next step is to 
	introduce a further modification to  $\psi_{1}$ to eliminate the two terms
	$$
	- {R^2 \over 2(h^2+ R^2)^2 } \left( {2h^2\over h^2+ R^2} +1    \right) \Gamma_{\ve \mu} \quad {\mbox {and}} \quad  {4 R^3 \over h(h^2+ R^2)^{3\over 2}} {{\mbox {Re} } (z^3 ) \over (\ve^2 \mu^2 + |z|^2)^2}.
	$$
	For the first one, we observe that
	\begin{align*}
		\Delta (c_2 |z|^2 \Gamma_{\ve \mu} ) &  - {R^2 \over 2(h^2+ R^2)^2 } \left({2 h^2\over h^2+ R^2} +1  \right) \Gamma_{\ve \mu} \\
		&=\left( 4 c_2  - {R^2 \over 2(h^2+ R^2)^2 } \left({2h^2\over h^2+ R^2} +1  \right) \right)\Gamma_{\ve \mu} \\
		&+ 2 c_2 z \cdot \nabla \Gamma_{\ve \mu} + c_2 |z|^2 \Delta \Gamma_{\ve \mu}
	\end{align*}
	and choose $c_2$ as
	$$
	c_2  = {R^2 \over 8 (h^2+ R^2)^2 } \left({2 h^2\over h^2+ R^2} +1  \right).
	$$
	To correct the second term , we introduce
	$$
	h_1 (s ) = s^3 \int_s^1 {dx \over x^7} \int_0^x {\eta^7 \over (\ve^2 \mu^2 + \eta^2)^2}  \, d \eta.
	$$
	It solves
	$$
	h_1'' +{1\over s} h_1' - {9 \over s^2} h_1 + {s^3 \over (\ve^2 + s^2)^2 }=0,
	$$
	it is smooth and uniformly bounded as $\ve \to 0$, and $h_1 (s ) = O(s^3)$, as $s \to 0$.
	Writing  $z = |z| e^{i\theta}$, we have that
	\be\label{H1}
	H_{1}  (z ) := h_1 (|z|) \cos 3 \theta \quad {\mbox {solves}} \quad \Delta_z \left( H_{1}   \right) + {{\mbox {Re} } (z^3) \over (\ve^2 \mu^2 + |z|^2 )^2} = 0.
	\ee
	For any $\mu >0$, we define the approximate regularization for \eqref{g2} as follows
	\be \label{defpsi2}\begin{aligned}
		\psi_{\mu} (z) &= \Gamma_{\ve \mu} (z) \left( 1+c_1 z_1 + c_2 |z|^2\right) +   {4 R^3 \over h(h^2+ R^2)^{3\over 2}}  H_1 (z) ,\\
		c_1 & ={1\over 2}  {R h \over (h^2+R^2)^{3\over 2} }\\
		c_2 &= {R^2 \over 8 (h^2+ R^2)^2 } \left({2 h^2\over h^2+ R^2} +1  \right)
	\end{aligned}
	\ee
	and $H_1$ as in \eqref{H1}. 
	With this definition, we have that
	\begin{equation}\label{ee1}
		L (\psi_{\mu} ) (z)= \Delta \Gamma_{\ve \mu}  + {4 R (3h^2+ R^2) \over h(h^2+ R^2)^{3\over 2} }  { \ve^2 \mu^2 \, z_1 \over (\ve^2 \mu^2 + |z|^2 )^2}  
		+ E_3(z)  
	\end{equation}
	where $E_3$  is a
	smooth function, which is uniformly bounded as $\ve \mu \to 0$.

	In the original variables $x$, the function $\psi_\mu$ in \eqref{defpsi2} reads as
	\be \label{psi2}
	\Psi_{ \mu , P} (x) = \psi_\mu \left(A[P]^{-1} (x-P) \right)
	\ee
	where $A[P]$ is the matrix introduced in the change of variables \eqref{matrixdef}. 
	The function $\Psi_{P, \mu} (x) $ is smooth and represents a good 
	approximate Green's function for the operator $\nabla \cdot (K \nabla \cdot )$ in $\R^2$. We will use it as building blocks for the construction of a solution to \eqref{PB1}.

	\subsection{Approximate stream function for $N$-helical filaments.} \ \ 
	Let $N$ be a fixed integer and consider $N$ points $P_1, \ldots , P_N$ of the form \eqref{points} and satisfying \eqref{points1}. For any $j=1, \ldots , N$, we write
	$$
	P_j = (a_j , b_j),
	$$
	we  fix positive constants $\mu_j$ and  define $$
	\Psi_j (x) =\Psi_{P_j, \mu_j} (x)
	$$
	where $\psi_{P_j , \mu_j}$ is defined in \eqref{psi2}. Hence, for each $j$ we have
	\be \label{Lpsij}
	L (\Psi_j ) = \Delta \Gamma_{\ve \mu_j}  + {4 R_j (3h^2+ R_j^2) \over h(h^2+ R_j^2)^{3\over 2} }  { \ve^2 \mu_j^2 \, z_1 \over (\ve^2 \mu_j^2 + |z|^2 )^2}  
	+ E_{3,j}(z)  
	\ee
	where
	$$
	R_j= \sqrt{a_j^2 + b_j^2}, \quad z = A[P_j]^{-1} (x-P_j).$$
	The stream function of $N$-helical filaments looks at main order as a superposition of stream functions $\Psi_j$ associated to each helical filament. Since the relative distance of the points $P_j$ is of order $|\log \ve |^{-1}$, we multiply each $\Psi_j$ by a cut-off function to get 
	$$
	\eta_0 (x) \sum_{j=1}^N \, \kappa_j \,  \,  \Psi_j (x) 
	$$
	where 
	\begin{equation}\label{defeta}
		\eta_{0} (x ) = \eta \left( |x-(r_0 , 0)|  \right),
	\end{equation}
	with $\eta $ a fixed smooth function with
	\begin{equation}\label{defeta2}
		\eta(s) = 1 , \quad {\mbox {for}} \quad s \leq {1 \over 2}, \quad \eta(s) = 0  , \quad {\mbox {for}} \quad s \geq 1.
	\end{equation}

	Using the notation introduced in \eqref{Lpsij}, we have
	\begin{align*}
		L\left( \eta_0 \sum_{j=1}^N \,  \, \kappa_j \, \Psi_j  \right)&=
		\eta_0 \sum_{j=1}^N \, \kappa_j \,  \,   \left( \Delta \Gamma_{\ve \mu_j}  + {4 R_j (3h^2+ R_j^2) \over h(h^2+ R_j^2)^{3\over 2} }  { \ve^2 \mu_j^2 \, z_1 \over (\ve^2 \mu_j^2 + |z|^2 )^2}  \right) \\
		&+g(x) , \quad {\mbox {where}} \\
		g(x)&= \eta_0 \sum_{j=1}^N \, \kappa_j \, \, E_{3,j}(z)  +
		\sum_{j=1}^N \kappa_j \left[  \, L( \eta_{0} \Psi_j ) - \eta_0 L(\Psi_j)   \right]
	\end{align*}
	The function $g$ has compact support and satisfies
	$$
	\| g (x) \|_{L^\infty (\R^2 ) } \leq C_{\delta_1}
	$$
	for some positive constant which depends on $\delta_1$.

	\medskip
	It is convenient to slightly modify the ansatz $\eta_{0} \sum_{j=1}^N \, \kappa_j \,  \,  \Psi_j $ adding a term which is defined globally in the entire space $\R^2$ to cancel $g(x)$. 
	Let $H_{2\ve} (x)$ solve
	\be \label{H2e}
	L (H_{2\ve}) + g = 0 , \quad {\mbox {in}} \quad \R^2.
	\ee
	For a smooth function $h(x)$ satisfying the decay condition
	$$
	\| h\|_\nu\, :=\, \sup_{x\in \R^2} (1+|x|)^\nu|h(x)|\, <\,+\infty \, ,
	$$
	for some $\nu >2$, there exists a solution $\psi(x)$ to problem 
	\be \label{outer}
	L(\psi) + h = 0 , \quad {\mbox {in}} \quad \R^2
	\ee
	which is of class $C^{1,\beta}(\R^2)$ for any $0<\beta<1$,  and defines a linear operator $\psi = {\mathcal T}^o (g) $ of $g$
	and satisfies the bound
	\be\label{estimate}
	|\psi(x)| \,\le \,     C{ \| h\|_\nu} (1+ |x|^2),
	\ee
	for some positive constant $C$. The proof of this fact can be found in Proposition \ref{prop2}. 
	
	\medskip
	Using this result we obtain that the solution $H_{2\ve}$ to \eqref{H2e} satisfies the estimate
	\begin{equation*}
		| H_{2\ve } (x) | \leq C_{\delta_1} (1+ |x|^2).
	\end{equation*}
	Besides, observe that such  solution is given up to the addition of a constant.
	We define the  function $H_{2\ve} (x)$  to be the one which furthermore satisfies
	$$
	H_{2\ve} ( (x_0, 0) ) = 0.
	$$
	With this is mind, we get to the definition of a first approximate stream function for $N$-helical filaments
	\be \label{psi0}
	\Psi_0 (x) = \eta_0 (x) \sum_{j=1}^N \, \kappa_j \,  \, \Psi_j (x) +H_{2\ve}(x),
	\ee
	so that
	\begin{equation} \label{psi01}
		L\left( \Psi_0   \right)= \eta_0
		\sum_{j=1}^N \, \kappa_j \, \,   \left( \Delta \Gamma_{\ve \mu_j}  + {4 R_j (3h^2+ R_j^2) \over h(h^2+ R_j^2)^{3\over 2} }  { \ve^2 \mu_j^2 \, z_1 \over (\ve^2 \mu_j^2 + |z|^2 )^2}  \right) .
	\end{equation}
	We recall that the definition of $\eta_0$ is given in \eqref{defeta}.
	The function $\Psi_0$ is defined in the whole $\R^2$ and it is smooth. Recalling that $\Psi_j (x) =\Psi_{P_j, \mu_j} (x)$, we observe that its definition depends on certain parameters: the points $P_1, \ldots , P_N$ and the scaling positive parameters $\mu_1, \ldots , \mu_N$. We now proceed to define the non-linearity $F$ in \eqref{PB1} and the first approximate solution to \eqref{PB1}. More specifically, we will define the scaling parameters $\mu_i$ as functions of the points $P_i$, assuming $P_i$ have the form \eqref{points}-\eqref{points1}.

	\section{Choice of the non-linearity and construction of the first approximate solution} \label{sec4}

	In this section we define a  nonlinearity $F$ in \eqref{PB1} with the property that the vorticity $W$, defined as
	$$
	W(x) = F (\Psi - {\alpha \over 2} |\log \ve | |x|^2 )$$
	satisfies \eqref{vort}, namely
	$$
	W(x) \sim 8 \pi \sum_{j=1}^N \kappa_j \delta_{P_j}, \quad {\mbox {as}} \quad \ve \to 0.
	$$
	To this purpose we first identify the form of
	the function
	$$
	\Psi_0 - {\alpha \over 2} |\log \ve | |x|^2
	$$
	near each point $P_j$, and we choose the scaling parameters $\mu_j$ in terms of the points $P_1, \ldots , P_N$. It will turn out that a convenient choice for $\mu_j$ gives their size of the order
	$$
	\log \mu_j \sim \log |\log \ve| , \quad \ass \ve \to 0.$$
	
	\subsection{Choice of $\mu_j$ in the definition of $\Psi_0$.}
	
	Fix $i \in \{ 1, \ldots , N\}$ and let $A_i = A[P_i]$ the matrix defined in \eqref{matrixdef}. Under our assumptions on the points $P_i$ in \eqref{points}-\eqref{points1}, and recalling that
	$$
	P_i = (r_0 + s, 0) + {1\over |\log \ve |} \hat P_i , \quad P_i = (a_i , b_i ) , \quad \hat P_i = (\hat a_i , \hat b_i ),
	$$
	we have
	\begin{align*}
		A_i^{-1} &= \left( \begin{matrix} {a_i \sqrt{h^2 + R_i^2} \over R_i h} & {b_i \sqrt{h^2 + R_i^2} \over R_i h}\\
			-{b_i \over R_i} & {a_i \over R_i}\end{matrix}\right) =\left( \begin{matrix} {\sqrt{h^2 + r_0^2} \over h} & 0 \\ 0 & 1 \end{matrix} \right) + {\log |\log \ve | \over |\log \ve |} \tilde A_i,\\
		{\mbox {and for } } i \not= j  & \\
		A_j^{-1}A_i&=
		\left(\begin{matrix} 
			1 &  0  \\
			0 & 1
		\end{matrix} \right) +\frac {\log |\log \ve |}{|\log \ve|} I_{ij}, 
	\end{align*}
	where $\tilde A_i$ and $I_{ij}$ are $2 \times 2$ matrices whose entrances are smooth functions of $(s, \hat a_i , \hat b_i)$, and $(s, \hat a_i , \hat b_i, \hat a_j , \hat b_j)$ respectively, which are uniformly bounded as $\ve \to 0$.

	\medskip{}
	We use these matrices to introduce useful changes of variables around each point $P_i$. 
	Take $\delta >0$ to be a fixed positive number and consider the inner region around $P_i$ to be given by 
	\be \label{innreg1}
	|A_i^{-1}(x-P_i)| <{\delta \over {|\log \ve|}}.
	\ee
	We take $\delta \leq {\sqrt{h^2 + r_0^2} \over h} \, {d \over 4}$, where $d $ was fixed in \eqref{defd}, so that for $x \in \R^2 $ satisfying  \eqref{innreg1} then $\eta_0 (x) =1$. 
	
	To represent a point in this region it is  convenient to use the change of variables  
	\be\label{varinner}
	x-P_i = A_i z, \quad z = \ve \mu_i y.
	\ee
	Hence
	$$
	|z| < {\delta \over |\log \ve |}.
	$$
	In this region we have the following expansion for $\Psi_0$, as $\ve \to 0$,
	\be\label{uno}
	\begin{aligned}
		&\Psi_0 (x )-\frac{\alpha}2 |\log \ve | |x|^2= \kappa_i \Gamma_0 (y) -\frac{\alpha}2 |\log \ve | |P_i|^2 - 4 \kappa_i \log \ve \mu_i \\
		&- 4 \kappa_i\log \ve \mu_i \left( c_{1,i} y_1 \, \ve \mu_i + c_{2,i} \ve^2 \mu_i^2 |y|^2 \right) + \kappa_i \frac{4R_i^3}{h (h^2+R_i^2)^{3/2}}H_{1i} (\ve \mu_i y)\\
		&-\alpha |\log \ve | {R_i h\over \sqrt{h^2+ R_i^2} }\ve \mu_i y_1 -\frac{\alpha}2 |\log \ve |  \ve^2 \mu_i^2 |A_i y|^2\\
		& +\kappa_i (c_{1,i} \ve \mu_i y_1+c_{2,i} \ve^2 \mu_i^2|y|^2)\Gamma_0(y)+ H_{2\ve} (\ve \mu_i A_iy+P_i)\\
		& + \sum_{j\not= i} \kappa_j \log {8\over |A_j^{-1} (P_i-P_j)|^4} \left( 1 + c_{1,j}[ A_j^{-1} (P_i - P_j)]_1 + c_{2,j} |A_j^{-1} (P_i - P_j) |^2 \right) \\
		&- \sum_{j\not= i} \kappa_j 4 {A_j^{-1} (P_i - P_j) \cdot A_j^{-1} A_i y \over |A_j^{-1} (P_i - P_j)|^2} \ve \mu_i \left( 1 + c_{1,j} [A_j^{-1} (P_i - P_j)]_1 + c_{2,j} |A_j^{-1} (P_i - P_j) |^2 \right)\\
		&+ \sum_{j\not=i} \kappa_j 
		\log {8\over |A_j^{-1} (P_i-P_j)|^4}  \{c_{1,j} [A_i^{-1} A_j y]_1 +2c_{2,j}A_j^{-1} (P_i-P_j)\cdot  A_j^{-1} A_i y\} \, \ve \mu_i\\
		&+ O\left( \log(|A_j^{-1} (P_i-P_j)|)|A_j^{-1} A_i y|^2 \ve^2\mu_i^2
		\right) 
	\end{aligned}
	\ee
	where
	the constants $c_{1,j}$, $c_{2,j}$ are defined as in \eqref{defpsi2}, namely 
	\begin{equation}\label{defc1j}
		\begin{aligned}
			c_{1,j} & ={1\over 2}  {R_j h \over (h^2+R_j^2)^{3\over 2} }\\
			c_{2,j} &= {R_j^2 \over 8 (h^2+ R_j^2)^2 } \left({2 h^2\over h^2+ R_j^2} +1  \right),
		\end{aligned}
	\end{equation}
	with $R_j= \sqrt{a_j^2 + b_j^2}$.

	To get expansion \eqref{uno} we have used that for small $z$, $|z| <{\delta \over |\log \ve |}$, one has
	\begin{align*}
		&\Gamma_{\ve\mu_j} (A_j^{-1}[A_i z-(P_j-P_i)] ) =  \log \frac{8 }{(\ve^2 \mu_j^2 + |A_j^{-1}[A_i z-(P_j-P_i)] |^2)^2} \\
		&	= \log \frac{8 }{|A_j^{-1}(P_j-P_i)|^4} - 2 \log \left( 1-2{A_j^{-1}(P_j-P_i) \cdot A_j^{-1}A_i z\over |A_j^{-1}(P_j-P_i)|^2} + {|A_j^{-1}A_iz|^2 + \ve^2 \mu_j^2 \over |A_j^{-1}(P_j-P_i)|^2} \right) \\
		&=	\log \frac{8}{|A_j^{-1}(P_j-P_i)|^4} + 4{ A_j^{-1}(P_j-P_i)   \over |A_j^{-1}(P_j-P_i)|^2} \, \cdot A_j^{-1}A_i z + O \left({\ve^2 \mu_j^2  + |A_j^{-1}A_i z|^2  \over |A_j^{-1}(P_j-P_i)|^2}\right)
	\end{align*} 
	as $\ve \to 0$.
	
	Observe that in the region  we are considering, we have the validity of the following expansion
	\begin{align*}
		H_{2\ve}(P_i) &=H_{2\ve}(P_i) +\ve\mu_i(A_iy)\cdot\nabla H_{2\ve}(P_i) +O(\ve^2\mu_i^2|y|^2)\\
		H_{1i}(\ve\mu_iy)&=O(\ve^3\mu_i^3|y|^3), \quad {\mbox {as}} \quad \ve \to 0.
	\end{align*}
	We now define the scaling parameters $\mu_i$, to eliminate part of the zero-mode term of the expression in \eqref{uno}. More precisely we take $\mu_i$ to be given by the relations
	\be\label{mu1}
	\begin{aligned}
		2\kappa_i \log \mu_i =& \sum_{j\not= i} \kappa_j \log {8\over |A_j^{-1} (P_i-P_j)|^4} \Bigl( 1 + c_{1,j} [A_j^{-1} (P_i - P_j)]_1  \\
		&+c_{2,j} |A_j^{-1} (P_i - P_j) |^2 \Bigl) +H_{2\ve}(P_i).
	\end{aligned}
	\ee
	Since the points $P=(P_1, \ldots , P_N)$ satisfy \eqref{points}-\eqref{points1} we recognize that 
	$$\log\mu_i^2 = \log(|\log \ve|) \, m_i (P), \quad {\mbox {as}} \quad \ve \to 0
	$$
	where $m_i(P)$ are smooth functions, which are uniformly bounded together with their derivatives, as $\ve \to 0$.
	We define
	\be \label{defmu}
	\mu =\max_{i=1, \ldots , N} \, \mu_i.
	\ee
	We insert \eqref{mu1} into \eqref{uno} and 
	we eventually obtain
	\begin{align*}
		&{1 \over \kappa_i} \Bigl( \Psi_0 (x )-\frac{\alpha}2 |\log \ve | |x|^2 \Bigl)= (1+ c_{1,i} \, \ve \mu_i \, y_1 + c_{2,i} \ve^2 \mu_i^2 |y|^2) \, \Gamma_0 (y) \\
		&- \frac{\alpha}{2\kappa_i} |\log \ve | |P_i|^2 -   4\log \ve - 2   \log \mu_i  \\
		&+ \ve y_1 \mu_i\left[ |\log \ve | \,   \left(  4   c_{1,i} \, 
		-\alpha  {hR_i \over \kappa_i\sqrt{h^2+ R_i^2} } \right) \right. - 4c_{1,i} \log\mu_i\\
		&- \sum_{j\not= i} {\kappa_j \over \kappa_i} 4 {[A_j^{-1} (P_i - P_j) \cdot ]_1  \over |A_j^{-1} (P_i - P_j)|^2}  \left( 1 +  c_{1,j} [A_j^{-1} (P_i - P_j)]_1 + c_{2,j} |A_j^{-1} (P_i - P_j) |^2 \right)\,  \\
		& \sum_{j\not=i} {\kappa_j \over \kappa_i}
		\log {8\over |A_j^{-1} (P_i-P_j)|^4}  \{c_{1,j}  +2c_{2,j}[A_j^{-1} (P_i-P_j)]_1\} \\
		&\left. +\frac 1{\kappa_i}(A_i (1,0)^T)\cdot\nabla H_{2\ve}(P_i)\right]\\
		&+ \ve y_2 \mu_i \left[  
		- \sum_{j\not= i} {\kappa_j \over \kappa_i} 4 {[ A_j^{-1} (P_i - P_j) \cdot ]_2  \over |A_j^{-1} (P_i - P_j)|^2}  \left( 1 +  c_{1,j} [A_j^{-1} (P_i - P_j)]_1 + c_{2,j} |A_j^{-1} (P_i - P_j) |^2 \right)\, \right.\\
		& +\left. \sum_{j\not=i} 2{\kappa_j \over \kappa_i}
		\log {8\over |A_j^{-1} (P_i-P_j)|^4}  c_{2,j}[A_j^{-1} (P_i-P_j)]_2 +\frac 1{\kappa_i}(A_i (0,1)^T)\cdot\nabla H_{2\ve}(P_i)\right]\\
		&+ 
		O\left( \log(|A_j^{-1} (P_i-P_j)|)|y|^2 \ve^2\mu_i^2
		\right) 
	\end{align*}
	
	\subsection{Choice of the non-linearity $F$}
	
	We now choose the nonlinearity $F$ in \eqref{PB1} which gives a vorticity $W$ satisfying \eqref{vort}.  We let $F$
	\be \label{defF}
	F(\Psi -\frac{\alpha}2 |\log \ve| |x|^2 ) = \sum_{j=1}^N \ve^{2-\frac{\alpha}{2\kappa_j}R_j^2}  \kappa_j F_j \left( \frac 1{\kappa_j}(\Psi -\frac{\alpha}2 |\log \ve| |x|^2) \right)
	\ee
	where
	\be \label{defF1}
	F_j (s) = \eta^j ( s)f(s),\quad \mbox{where}\quad f(s)= e^{ s}.
	\ee
	Here $\eta^j$ are smooth cut-off functions defined as follows.

	\medskip
	Consider the boundary of the inner region around $P_i$, as defined in \eqref{innreg1}. Using the variable $y$ in \eqref{varinner}, this boundary is defined by $|y|=\delta/({\mu_i\ve|\log\ve|})$. On this boundary we have
	\begin{align*}
		{ 1 \over \kappa_i} &\Bigl( \Psi_0 (x )-\frac{\alpha}2 |\log \ve | |x|^2 \Bigl)= \\
		&-\frac{\alpha R_i^2}{2\kappa_i} |\log \ve |+4 \log|\log \ve|+2\log\mu_i+\log 8-4\log\delta+o(1)
	\end{align*}
	where $o(1)$ is respect to $\ve\to 0$.
	Then we choose the cutoff function $\eta^i$ such that 
	\begin{equation}\label{cutoff}
		\begin{aligned}
			\eta^i(s)&=1, \quad {\mbox {for}} \quad  s\geq -\frac{\alpha R_i^2}{2\kappa_i} |\log \ve |+4 \log|\log \ve|+2\log\mu_i +\log8+2d_{i,\ve} \\
			\eta^i(s)&=0 , \quad {\mbox {for}} \quad s\leq -\frac{\alpha R_i^2}{2\kappa_i} |\log \ve |+4 \log|\log \ve|+2\log\mu_i +\log8 +d_{i,\ve}
		\end{aligned}
	\end{equation}
	for suitable $d_{i,\ve}=-4\log \delta+o(1)$ so that 
	$$
	\eta^i\left({ 1 \over \kappa_i} \Bigl( \Psi_0 (x )-\frac{\alpha}2 |\log \ve | |x|^2 \Bigl)\right)=1\quad \mbox{for}\quad |A_i^{-1}(x-P_i)|\leq \frac{\delta^2}{|\log\ve|}
	$$
	and 
	$$
	\eta^i\left({ 1 \over \kappa_i} \Bigl( \Psi_0 (x )-\frac{\alpha}2 |\log \ve | |x|^2 \Bigl)\right)=0\quad \mbox{for}\quad |A_i^{-1}(x-P_i)|\geq \frac{\delta}{|\log\ve|}.
	$$
	Here $\delta$ is independent of $\ve$ as in \eqref{innreg1} and can be taken smaller if needed.
	
	
	\section{Estimate of the error function}\label{sec5}
	Let us define  the error function to be
	\be \label{error}
	S[\Psi ] (x) = L (\Psi) + \sum_{j=1}^N \ve^{2-\frac{\alpha}{2\kappa_j}R_j^2}   \kappa_j F_j \left( \frac 1{\kappa_j}(\Psi -\frac{\alpha}2 |\log \ve| |x|^2 )\right).
	\ee
	A solution to \eqref{PB1} would correspond to a smooth function $\Psi$ such that
	$$
	S[\Psi ] (x) = 0 \quad x\in \R^2.
	$$
	Our purpose is to estimate 
	$$
	S[\Psi_0 ] (x) \quad {\mbox {for}} \quad  x \in \R^2
	$$
	where $\Psi_0$ is the approximate stream function for the $N$-helical filaments introduced in \eqref{psi0}. 
	
	In order to do so, we shall first analyze $S[\Psi_0]$ 
	in regions that are close to each vortex point $P_j$, and then  in the region which is far from all the points $P_1, \ldots , P_N$. Let us be more precise.

	\medskip
	We split the inner region around $P_i$ as described in \eqref{innreg1}, into two parts
	$$ |A_i^{-1}(x-P_i)|\leq \frac{\delta^2}{|\log\ve|} \quad {\mbox {and}} \quad \frac{\delta^2}{|\log\ve|}\leq |A_i^{-1}(x-P_i)|\leq \frac{\delta}{|\log\ve|}.$$
	Assume first that $|A_i^{-1}(x-P_i)|\leq \frac{\delta^2}{|\log\ve|} $. According to \eqref{cutoff}, we have $\eta^i =1 $
	and $\eta^j =0$ for $j\not= i$, 
	so that the non-linear term in the expression of $S(\Psi_0)$ becomes
	\be \label{rhsf}
	\begin{aligned}
		&\sum_{j=1}^N \ve^{2-\frac{\alpha}{2\kappa_j}R_j^2}   \kappa_j F_j \left( \frac 1{\kappa_j}(\Psi_0-\frac{\alpha}2 |\log \ve| |x|^2 )\right)\\
		=  &\ve^{2-\frac{\alpha}{2\kappa_i}R_i^2}\,  \kappa_i f\left( { 1 \over \kappa_i} \Bigl( \Psi_0 (x )-\frac{\alpha}2 |\log \ve | |x|^2 \Bigl)\right)\\
		&= \frac{\kappa_i}{\ve^2\mu_i^2}U (y)  e^{[ c_{1,i} \ve \mu_i y_1 + c_{2,i} \ve^2 \mu_i^2 |y|^2] \, \Gamma_0 (y) } \, e^{ \ve [ {\mathcal A}_{1,i} (P) y_1 \mu_i + {\mathcal A}_{2,i} (P)  y_2 \mu_i ] } \\
		&\times \exp\left[   O\left( \log(|\log \ve |)|y|^2 \ve^2\mu_i^2
		\right)
		\right]\\
		&= \frac{\kappa_i}{\ve^2\mu_i^2}U (y)  \Biggl[ 1+ \ve \mu_i y_1 \left(c_{1,i} \Gamma_0 (y) + {\mathcal A}_{1,i} (P) \right) + \ve \mu_i y_2  {\mathcal A}_{2,i} (P)  \\
		& \quad  \quad \quad  \quad \quad  + \ve^2 \mu_i^2 c_{2,i} |y|^2 \Gamma_0 (y) +O\left( \log(|\log \ve |)|y|^2 \ve^2\mu_i^2
		\right) \Biggl]
		&   \\
		&{\mbox {with}}\\
		&{\mathcal A}_{1,i} (P): =|\log \ve | \,  \left(  4   c_{1,i} \, 
		-\alpha  {h R_i \over \kappa_i \sqrt{h^2+ R_i^2} } \right) - 4c_{1,i} \log\mu_i \\
		&- 4\sum_{j\not= i} {\kappa_j \over \kappa_i} {[A_j^{-1} (P_i - P_j)]_1  \over |A_j^{-1} (P_i - P_j)|^2}  
		+ \sum_{j\not=i} {\kappa_j \over \kappa_i} \log {8\over |A_j^{-1} (P_i - P_j)|^4} c_{1,j} \\
		&+{\bf Y}_1(P) \\
		&{\mathcal A}_{2,i} (P):=  - \sum_{j\not= i} {\kappa_j \over \kappa_i} 4 {[A_j^{-1} (P_i - P_j)]_2  \over |(A_j^{-1} (P_i - P_j))|^2} 
		+{\bf Y}_2(P)
	\end{aligned}
	\ee
	where ${\bf Y}_1(P)$ and ${\bf Y}_2(P)$ are smooth functions,  uniformly bounded  as $\ve \to 0$  for points  $P= (P_1, \ldots , P_N) $   satisfying  \eqref{points}-\eqref{points1}. 
	
	\medskip
	In the expression of ${\mathcal A}_{1,i}$ the term $ - 4c_{1,i} \log\mu_i +\sum_{j\not=i} {\kappa_j \over \kappa_i} \log {8\over |A_j^{-1} (P_i - P_j)|^4} c_{1,j} $ is a smooth function of the points $P_j$ which can be described as
	$\log |\log \ve | {\bf Y}_1 (P)$, where
	${\bf Y}_1 (P)$ denotes again a smooth function,  uniformly bounded  as $\ve \to 0$  for points  $P= (P_1, \ldots , P_N) $   satisfying  \eqref{points}-\eqref{points1}. 
	
	Besides, if we insert the definition of   $c_{1,i} $ as given in \eqref{defc1j}, we write $\mathcal A_{1,i} (P)$  as
	\be\label{aa1}
	\begin{aligned}
		{\mathcal A}_{1,i} (P) &= |\log \ve | \,   \left(  2  {h R_i \over  \sqrt{(h^2+ R_i^2)^3} }   \, 
		-\alpha  {h R_i \over \kappa_i\sqrt{h^2+ R_i^2} } \right) 
		\\
		&- \sum_{j\not= i} {\kappa_j \over \kappa_i} 4 { [A_j^{-1} (P_i - P_j)]_1  \over |A_j^{-1} (P_i - P_j) |^2} + \log|\log\ve| \, {\bf Y}_1 (P)
	\end{aligned}
	\ee
	where again ${\bf Y}_1(P)$ denotes an explicit smooth function,  uniformly bounded  as $\ve \to 0$  for points  $P= (P_1, \ldots , P_N) $   satisfying  \eqref{points}-\eqref{points1}. 
	
	\medskip
	For later purpose it is relevant to observe that
	\be \label{aa2}
	{\mathcal A}_{1,i} = \log |\log \ve | \, {\bf \overline{Y}}_1 (P), \quad {\mathcal A}_{2,i} = |\log \ve | \, {\bf \overline{Y}}_2 (P)
	\ee
	under the assumption that the points $P_i$ satisfy \eqref{points}-\eqref{points1}. As before ${\bf \overline{Y}}_1(P)$ and ${\bf \overline{Y}}_2 (P)$ denote an explicit smooth functions,  uniformly bounded  as $\ve \to 0$  for points  $P= (P_1, \ldots , P_N) $   satisfying  \eqref{points}-\eqref{points1}.

	\medskip
	A direct computation shows that in the region  $\frac{\delta^2}{|\log\ve|}\leq |A_i^{-1}(x-P_i)|\leq \frac{\delta}{|\log\ve|}$ we have
	\begin{align*}
		\sum_{j=1}^N &\ve^{2-\frac{\alpha}{2\kappa_j}R_j^2}   \kappa_j F_j \left( \frac 1{\kappa_j}(\Psi_0-\frac{\alpha}2 |\log \ve| |x|^2 )\right)\\
		&= \ve^{2-\frac{\alpha}{2\kappa_i}R_i^2}f\left( { 1 \over \kappa_i} \Bigl( \Psi_0 (x )-\frac{\alpha}2 |\log \ve | |x|^2 \Bigl)\right)\\
		&=O(\ve^2\mu_i^2|\log\ve|^4 ). 
	\end{align*}
	On the other hand we recall from \eqref{psi01} that 
	\be \label{ee0}
	L(\Psi_0)= \eta_0 \sum_{j=1}^N \, \kappa_j \, \,   \left( \Delta \Gamma_{\ve \mu_j}  + {4 R_j (3h^2+ R_j^2) \over h(h^2+ R_j^2)^{3\over 2} }  { \ve^2 \mu_j^2 \, z_1 \over (\ve^2 \mu_j^2 + |z|^2 )^2}  \right)
	\ee
	in the variable $z=(x-P_i)/\ve\mu_i$
	$$
	\begin{aligned}
		&L(\Psi_0)= \eta_0  \kappa_i \, \,   \left( \Delta \Gamma_{\ve \mu_i}  + {4 R_i (3h^2+ R_i^2) \over h(h^2+ R_i^2)^{3\over 2} }  { \ve^2 \mu_j^2 \, z_1 \over (\ve^2 \mu_i^2 + |z|^2 )^2}  \right) \\
		&+\eta_0\sum_{j\neq i}^N \, \kappa_j \, \,   \left( \frac{8\ve^2 \mu_j^2}{(\ve^2 \mu_j^2+|A_j^{-1}A_iz+A_j^{-1}(P_j-P_i)|^2 )^2} \right. \\
		&+ \left. {4 R_j (3h^2+ R_j^2) \over h(h^2+ R_i^2)^{3\over 2} }  { \ve^2 \mu_j^2 \, [A_j^{-1}A_iz+A_j^{-1}(P_j-P_i)]_1 \over (\ve^2 \mu_j^2 + |A_j^{-1}A_iz+A_j^{-1}(P_j-P_i)|^2 )^2}  \right)
	\end{aligned}
	$$
	Therefore for the \emph{inner part} $|A_i^{-1} (x-P_i) |<\delta/|\log \ve|$, we have 
	\be \label{ee1}
	\begin{aligned}
		L(\Psi_0)&
		=-\frac {\kappa_i}{\ve^2\mu_i^2}\left[U(y)-{4 R_i (3h^2+ R_i^2) \over h(h^2+ R_i^2)^{3\over 2} }  { \ve \mu_i \, y_1 \over (1 + |y|^2 )^2}\right]+O(\ve^2\mu^2|\log\ve|^4)
	\end{aligned}
	\ee
	where $y$  is the variable introduced in \eqref{varinner}
	$$
	x-P_i = \ve \mu_i A_i   y.
	$$

	\medskip
	Combining \eqref{rhsf} and \eqref{ee1}, we conclude that in the region
	$|A_i^{-1} (x-P_i) |<\delta^2/|\log \ve|$ we have
	\be \label{ee2}
	\begin{aligned}
		S(\Psi_0 ) &=\frac{\kappa_i}{\ve^2\mu_i^2}U (y)  \Biggl[ \ve \mu_i y_1 \left(c_{1,i} \Gamma_0 (y)+{ R_i (3h^2+ R_i^2) \over 2h(h^2+ R_i^2)^{3\over 2} }  + {\mathcal A}_{1,i} (P) \right) \\
		&+ \ve \mu_i y_2  {\mathcal A}_{2,i} (P)   + \ve^2 \mu_i^2 c_{2,i} |y|^2 \Gamma_0 (y) +O\left( \log(|\log \ve |)|y|^2 \ve^2\mu_i^2
		\right) \Biggl] \\
		&+O(\ve^2\mu^2|\log\ve|^4)
	\end{aligned}
	\ee
	with $\mu $ given by \eqref{defmu}, and $S$ by \eqref{error}.
	
	\medskip
	Consequently, 
	we get that in the region  $\frac{\delta^2}{|\log\ve|}\leq |A_i^{-1}(x-P_i)|\leq \frac{\delta}{|\log\ve|}$ we have
	$$
	S(\Psi_0 ) (x)
	= O(\ve^2\mu^2|\log\ve|^4)
	$$
	
	\medskip
	Let us consider now the region defined by
	$$|A_i^{-1}(x-P_i)|>{\delta \over |\log \ve|}, \quad \forall i.
	$$
	In this outer region, all cut-off functions $\eta^i$ are zero, see \eqref{cutoff}, hence
	$
	S[\Psi_0 ] (x) = L (\Psi_0).$
	A direct inspection of \eqref{ee0} gives
	$$
	\left|L\left(\Psi_0 \right)\right|\leq C\frac{\ve^2\mu^2}{1+|x|^\nu}
	$$
	for some $C$ independent of $\ve$, $\nu>2$ and $\mu$ is defined in \eqref{defmu}.

	\medskip
	Combining all previous results, we conclude that: in the region
	$$
	|A_i^{-1}(x-P_i)|>{\delta \over |\log \ve|}, \quad \forall i =1, \ldots , N
	$$
	we have
	\be
	\label{SPsi02}
	|S(\Psi_0)|\leq C\frac{\ve^{1+\sigma}}{1+|x|^{\nu}}.
	\ee
	for some constant $C$ independent of $\ve$, $\nu>2$ and $\sigma\in (0,1)$.

	For $i=1, \ldots , N$, in the region
	$$
	|A_i^{-1}(x-P_i)|<{\delta \over |\log \ve|}, 
	$$
	we have
	\be
	\label{SPsi01}
	\ve^2\mu_i^2 S(\Psi_0)\leq C  {  \ve \mu_i \log|\log\ve| \over  (1 + |y|^{2+a})} 
	\ee
	for some constant $C$ independent of $\ve$ and  $a\in (0,1)$. Estimate \eqref{SPsi01} uses \eqref{aa1}-\eqref{aa2}.

	\medskip
	Estimates \eqref{SPsi02} and \eqref{SPsi01} will be crucial to carry on the inner-outer gluing procedure leading to an exact solution of \eqref{PB1}. This is what we discuss next.


	\section{The inner-outer gluing system}

	This section describes  the inner-outer gluing scheme to find an actual solution to \eqref{PB1}.
	We look  for a solution $\Psi(x) $
	of the equation
	\be\label{ecuacion}
	S[\Psi] := L[\Psi]  + F( \Psi )=0 \inn \R^2
	\ee
	where
	$$
	F(\Psi)=\sum\limits_{i=1}^N\ve^{2-{\alpha \over 2 \kappa_i}R_i^2}\kappa_i \eta^i
	f\left( \frac 1{\kappa_i}(\Psi -\frac{\alpha}2 |\log \ve| |x|^2) \right), \quad f(u)= e^u.
	$$
	Consider the approximate solution $\Psi_0(x)$ in \eqref{psi0}. The function $\Psi_0$ is defined in terms of scaling parameters  $\mu_1 , \ldots , \mu_N$ given by formula  \eqref{mu1} and points $P_1, \ldots , P_N$ satisfying \eqref{points}-\eqref{points1}. We refer to  Sections \S \ref{sec3} and \S \ref{sec4} for the construction of $\Psi_0$.
	We look for a solution $\Psi$ to \eqref{ecuacion} of the form
	\be \label{ansat}
	\Psi (x) = \Psi_{ 0} (x) + \varphi(x).
	\ee
	where $\varphi$ is  "smaller" than  $\Psi_0$. 
	The inner-outer gluing procedure starts with  choosing $\varphi$ of the form
	\be \label{ff1}
	\varphi(x)=   \sum\limits_{i=1}^N\eta_i (x) \phi_i\left (y \right)  +  \psi (x) .
	\ee
	Here $$\eta_i=\eta\left({|\log\ve||A_i^{-1}(x-P_i)|\over {\delta_1}}\right)$$ for some $\delta_1<\delta^2$, with $\delta$ fixed in \eqref{innreg1} and  $\eta$  defined in \eqref{defeta2},  and $y$ denotes the scaling variable
	$$y=\frac{A_i^{-1}(x-P_i)}{\ve \mu_i}.$$
	In terms of $\varphi(x)$, and using  the decomposition \eqref{ff1}, problem \eqref{ecuacion} takes the form
	$$
	\begin{aligned} S( \Psi_0 +\vp) &= 0 \quad \inn \quad \R^2, \\
		S( \Psi_0 +\vp)&=\
		\sum\limits_{i=1}^N\eta_i\big [  L[\phi_i]  +  F'( \Psi_0 ) (\phi_i + \psi)  + S (\Psi_0 ) +  N_{0 } (\sum\limits_{i=1}^N\eta_i \phi_i + \psi ) \big]
		\\
		& +    L[\psi]   + (1-\sum\limits_{i=1}^N\eta_i )\left[   F'( \Psi_0 ) \psi  +   E_0 +  N_{0 } (\sum\limits_{i=1}^N\eta_i \phi_i + \psi ) \right]\\
		&+ \sum\limits_{i=1}^N \left( L[\eta_i \phi_i ] -\eta_i L[\phi_i]  \right)
	\end{aligned}
	$$
	where
	$$
	\begin{aligned}
		N_{0 } (\varphi ) &=  F (\Psi_0 +\varphi) -
		F  (\Psi_0  ) - F' (\Psi_0 ) \varphi , \quad {\mbox {with}} \\
		F'(\Psi_0)& = \sum\limits_{i=1}^N\ve^{2-{\alpha \over 2 \kappa_i}R_i^2} \eta^i \, 
		f'\left(\frac 1{\kappa_i}(\Psi_0 -\frac{\alpha}2 |\log \ve| |x|^2)\right).
	\end{aligned}
	$$
	
	Thus $\Psi$ given by \equ{ansat}-\equ{ff1} solves \equ{ecuacion} if  $(\phi,\psi):=(\phi_1,\ldots,\phi_N,\psi)$ satisfies  the system of equations
	\be
	\begin{aligned}
		L[\phi_i]  +  F'( \Psi_0 ) (\phi_i + \psi)  &+ S (\Psi_0 )+  N_{0} (\sum\limits_{i=1}^N\eta_i \phi_i + \psi ) \, =\, 0,  \\
		\quad  {\mbox {for}} &\quad |A_i^{-1}(x-P_i)|< \frac{2\delta_1}{|\log \ve|},
	\end{aligned}
	\label{in}\ee
	and
	\be\begin{aligned}
		&  L[\psi]   + (1-\sum\limits_{i=1}^N\eta_i )\left[   F'( \Psi_0 ) \psi  +   E_0 +  N_{0 } (\sum\limits_{i=1}^N\eta_i \phi_i + \psi ) \right]\\
		&+ \sum\limits_{i=1}^N \left( L[\eta_i \phi_i ] -\eta_i L[\phi_i]  \right) =\  0 \inn \R^2.\end{aligned} \label{out}  \ee
	We will refer to problem \eqref{in} as the  {\it inner problem} and to \eqref{out} as the {\it outer problem}.
	
	\medskip
	Let us write \equ{in} in terms of the variable $y=\frac {A_i^{-1}(x-P_i)}{\ve\mu_i}$. From \eqref{Lz} we have
	$$
	L[\phi_i]   =   \frac 1{\ve^2\mu_i^2} \big[ \,\Delta_y \phi_i + \bar B_i(y)[\phi_i]  \,\big ]
	$$
	where $\bar B_i(y)=\ve^2\mu_i^2B(\ve\mu_iy)$ and $B$ is the operator given by \eqref{B0}.
	Using the estimate \eqref{rhsf}, we get
	$$
	\ve^2\mu_i^2F'( \Psi_0 )=e^{\Gamma_0(y)}+b_i(y)\quad \mbox{with}\quad \Gamma_0 (y) = \log {8 \over (1+ |y|^2)^2}
	$$
	where
	\be \label{biest0}
	\begin{aligned}
		b_i(y)&=e^{\Gamma_0(y)}  \Biggl[\ve \mu_i y_1 \left(c_{1,i} \Gamma_0 (y) + {\mathcal A}_{1,i} (P) \right) + \ve \mu_i y_2  {\mathcal A}_{2,i} (P)  \\
		& \quad  \quad \quad  \quad   + \ve^2 \mu_i^2 c_{2,i} |y|^2 \Gamma_0 (y) +O\left( \log(|\log \ve |)|y|^2 \ve^2\mu_i^2
		\right) \Biggl].
	\end{aligned}
	\ee
	Consequently using \eqref{aa1}-\eqref{aa2}, we have 
	\be \label{biest}
	b_i(y)=O\left(
	\frac{\ve \mu_i \log|\log\ve|}{1+|y|^{2+a}}\right), 
	\ee
	so the term $\log|\log\ve|$ assumes that the points $P_i$ satisfy \eqref{points}-\eqref{points1}.
	
	By estimates \eqref{SPsi01}, 
	we obtain, in the region $|y|<{2\delta_1 \over \mu_i\ve|\log\ve|}$,
	\begin{equation}\label{ba1}
		\tilde E_i :=\ve^2\mu_i^2 S[\Psi_0 ] \ = \
		O\left( \frac{ \ve \mu_i }  {1+|y|^{2+a}}\log|\log\ve|\right) 
	\end{equation}
	for $a \in (0,1)$.
	Similarly, using estimate \eqref{biest} for $b_i$, we get the expansion
	\begin{equation}\label{house2}
		\mathcal N_i (\vp) := \ve^2\mu_i^2 N_{0} (\vp)  = \frac 1{\kappa_i} ( e^{\Gamma_0(y)}   +   b_i(y))  \vp^2.
	\end{equation}
	Then, multiplying  the inner problem \equ{in} by $\ve^2\mu_i^2$, we get
	\be\label{inner}
	\Delta_y \phi_i  +   e^{\Gamma_0} \phi =-   B_i[\phi_i]  -  H_i(\phi, \psi)    \inn B_R
	\ee
	where 
	$R= \frac {2\delta_1}{\ve\mu_i|\log\ve|} $ ,
	$$
	H_i(\phi, \psi)  = {\mathcal N}_i \left(\sum\limits_{i=1}^N\eta_i \phi_i + \psi \right) +   \ttt  E_i +  (e^{\Gamma_0} + b_i ) \psi
	$$
	and
	\begin{equation}\label{defB}
		B_i[\phi_i] = \bar B_i(y)[\phi_i] + b_i(y) \phi_i.
	\end{equation}
	Note that
	\be \label{By}
	\begin{aligned}
		\bar B_i(y)&=\left( -2{R_i h \over  (h^2+ R_i^2)^{3/2}} \, (\ve \mu_iy_1) + O(|\ve \mu_iy|^2 )\right) \pp_{y_1 y_1} + O(|\ve \mu_iy|^2 ) \pp_{y_2 y_2} \\
		&-\left( 2{R_i \over h\sqrt{h^2+ R_i^2} } \ve \mu_iy_2 + O(|\ve \mu_iy|^2) \right) \pp_{y_1 y_2}\\
		&
		- \left( {\ve \mu_i R_i \over h\sqrt{h^2+ R_i^2} } \left({2h^2\over h^2+ R_i^2} +1 \right)  + O((\ve \mu_i)^2|y|) \right) \pp_{y_1 } \\
		&- \left(\frac{(\ve \mu_i)^2y_2}{h^2+R_i^2} \left({2 h^2\over h^2+ R_i^2} +1 \right)+O((\ve \mu_i)^3|y|^2)\right) \pp_{y_2 }.
	\end{aligned}
	\ee
	
	The idea is to solve equation \eqref{inner}, coupled with the outer problem \equ{out} in such a way that  $\phi_i$ has the size of the error 
	$\ttt E_i$ 
	with two powers less of decay in $y$, say
	$$
	(1+|y|) |D_y\phi_i(y)|+  |\phi_i(y)| \le  \frac {C \ve \mu_i \log|\log\ve| } {1+ |y|^{a}}. 
	$$
	Recall that the basic linear operator $\Delta_y \phi + e^{\Gamma_0} \phi $ in \eqref{inner} has $3$ dimensional kernel generated by the bounded functions
	\be \label{defZj}
	Z_i (y) = {\pp \Gamma_0 \over \pp y_i } , \quad i=1,2, \quad Z_0 (y) = 2+ y\cdot \nabla \Gamma_0 (y) .
	\ee
	This fact suggests that solvability of \eqref{inner} within the expected topologies depends on whether the right hand side does have component in the directions spanned by the $Z_i$. Instead of solving directly \eqref{inner}, we will solve the  auxilliary projected problem instead problem, for $i=1\ldots N$,
	\begin{equation}\label{innerp}
		\Delta_y \phi_i + e^{\Gamma_0 } \phi_i + B_i(\phi_i)+H_i(\phi,\psi)=\sum_{j=1}^2c_{ij}e^{\Gamma_0(y)}Z_j\quad \mbox{in}\quad B_R,
	\end{equation}
	for some constants $c_{ij}$. We solve \eqref{innerp} coupled with the outer problem \eqref{out}, which can be written  as
	\be\label{outer}
	L[\psi]  \ +  \,   G(\psi, \phi)  \ =\  0 \inn \R^2
	\ee
	where
	\be
	G(\psi, \phi) =
	V(x)\psi \, +   \, N^0( \phi  ) \, +  \, E^0 (x)\,
	+ \ \sum\limits_{i=1}^N A_i[\phi_i],
	\label{GG}\ee
	with
	$$
	\begin{aligned}
		V(x)\ = &\ (1-\sum\limits_{i=1}^N\eta_i )  F'( \Psi_0 ) , \quad   N^0(\vp) \  = \  (1-\sum\limits_{i=1}^N\eta_i )  N_{0 } (\vp )\\
		E^0 (x)\ = &\  (1-\sum\limits_{i=1}^N\eta_i) S[\Psi_0 ],\quad   A_i[\phi_i]\ =\  L[\eta_i] \phi_i\, +\,  K_{\ell j} (x)\pp_{x_\ell} \eta_i\pp_{x_j}\phi_i,
	\end{aligned}$$
	where $K_{\ell j}$ are the coefficients of the matrix $K$ defining the differential operator $L$, see  \eqref{defL}.
	By \eqref{SPsi02}, the following bounds hold
	\begin{equation}\label{house1}
		\begin{aligned} |V(x)|\ \le\   O((\ve\mu)^{2}|\log\ve|^4)  &,\quad  
			|N^0(\varphi)|  \ \le \   O((\ve\mu)^{2}|\log\ve|^4   |\varphi|^2), \\ 
			|E^0 (x)|\ &\le \  O( \ve^{1+b}).
		\end{aligned}
	\end{equation}
	In order to found a solution of \eqref{ecuacion}, we will need to solve
	$$
	c_{ij}=c_{ij}[\phi,\psi]=0\quad \mbox{for}\quad i=1\ldots N,\quad j=1,2.
	$$
	This can be achieved choosing properly $s^*$ and $\hat P_1, \ldots , \hat P_N$ in the form of the points $P_j$ as given in \eqref{points}, under the bounds \eqref{points1}.
	In Section \S \ref{appe}  we will establish linear results that are the basic tools to solve system \equ{innerp}-\equ{outer}. Section \S \ref{sette} is devoted to solve \equ{innerp}-\equ{outer} by means of a
	fixed point scheme, and Section \S \ref{otto} to adjust the points to get $c_{ij}=0$ for all $i=1, \ldots , N$, $j=1,2$.

	\bigskip
	
	\section{Linear theories}\label{appe}

	This section collects two results. The first one regards the solvability of the {\it outer linear theory}, and the second the {\it inner linear theory}. They have been obtained in \cite{ddmw2}. For completeness we state them here and give a sketch of their proofs in the Appendix \S \ref{appe1}.

	\subsection{Outer linear theory}

	Consider the Poisson equation for the operator $L$
	\be
	\label{louter}
	L[\psi]  +  g(x) = 0  \inn \R^2 ,
	\ee
	for a bounded function $g$. Here $L$ is the differential operator in divergence form defined in \eqref{defL}.

	We take functions  $g(x)$ that satisfy the decay condition
	$$
	\| g\|_\nu\, :=\, \sup_{x\in \R^2} (1+|x|)^\nu|g(x)|\, <\,+\infty \, ,
	$$
	where $\nu >2$. 
	\begin{prop}\label{prop2} (Proposition 7.1 in \cite{ddmw2})
		There exists a solution $\psi(x)$ to problem $\equ{louter}$, which is of class $C^{1,\beta}(\R^2)$ for any $0<\beta<1$,  that defines a linear operator $\psi = {\mathcal T}^o (g) $ of $g$
		and satisfies the bound
		\be\label{estimate}
		|\psi(x)| \,\le \,     C{ \| g\|_\nu} (1+ |x|^2),
		\ee
		for some positive constant $C$.
	\end{prop}

	\subsection{Inner linear theory}
	In this section we consider the problem
	\be\label{00}
	\begin{aligned}
		\Delta \phi + e^{\Gamma_0(y)}\phi  +  h(y)= 0 \inn \R^2.  
	\end{aligned}
	\ee
	For numbers $m>2$, $0<\beta <1$ we consider the following norms
	\be\label{norma} \begin{aligned}
		\|h\|_{m}  = & \sup_{y\in \R^n} (1+|y|)^m|h(y)|, \\
		\|h\|_{m,\beta}  = &\|h\|_{m}  +   (1+|y|)^{m + \beta}[h]_{B_1(y),\beta} ,
	\end{aligned}
	\ee
	where we use the standard notation
	$$
	[h]_{A,\beta} = \sup_{z_1,z_2\in A}  \frac {|h(z_1) -h(z_2)| } {|z_1-z_2|^\beta},
	$$
	and $A$ is a subset of $\R^2$.
	We recall the definition of the functions $Z_i(y)$ in \eqref{defZj}
	$$ Z_i(y)  =  \pp_{y_i}  \Gamma_0(y) , \  i=1,2, \quad  Z_0(y) =    2+  y\cdot \nn \Gamma_0(y)  .     $$

	\begin{lemma}\label{lemat}
		(Lemma 6.1 in \cite{ddmw2}). Given $m>2$ and $0<\beta < 1$,
		there exists a $C>0$ and a solution $\phi =  \TT [ h]$ of problem $\equ{00}$ for each $h$ with $\|h\|_m <+\infty$ that defines a linear operator of $h$ and satisfies the estimate
		\be\label{cota} \begin{aligned}
			& (1+|y|) | \nn \phi (y)|  +  | \phi (y)| \\  &
			\,  \le  \,  C \big [ \,  \log (2+|y|) \,\big|\int_{\R^2} h Z_0\big|  +    (1+|y|) \sum_{j=1}^2 \big|\int_{\R^2} h Z_j\big| \\
			&\qquad +  (1+|y|)^{2-m} \|h\|_{m}   \,\big ]. \end{aligned}
		\ee
		In addition, if
		$\|h\|_{m,\beta} <+\infty$, we have
		\be\label{cotaa} \begin{aligned}
			& (1+|y|^{2+\beta})  [D^2_y \phi]_{B_1(y),\beta}  +(1+|y|^2)  |D^2_y \phi (y)| \\  &
			\,  \le  \,  C \big [ \,  \log (2+|y|) \,\big|\int_{\R^2} h Z_0\big|  +    (1+|y|) \sum_{j=1}^2 \big|\int_{\R^2} h Z_j\big|  \\
			& \qquad +  (1+|y|)^{2-m} \|h\|_{m,\beta}   \,\big ]. \end{aligned}
		\ee
		
	\end{lemma}

	\bigskip
	We consider now the problem for a fixed number $\delta>0$ and a sufficiently large $R>0$ we consider the equation
	\be
	\Delta \phi  +  e^{\Gamma_0} \phi  + B_i[\phi]  + h(y)  = \sum_{j=0}^2  c_{ij} e^{\Gamma_0} Z_j \inn B_R
	\label{eee} \ee
	For a function $h$ defined in $A\subset \R^2$ we denote by $\|h\|_{m,\beta,A}$ the numbers defined in \equ{norma}
	but with the sup  taken with elements in $A$ only, namely
	\begin{align*}
		\|h \|_{m, A} &=  \sup_{y\in A} | (1+|y|)^m h(y)|] ,\\     \|h \|_{m,\beta, A} = & \sup_{y\in A} (1+|y|)^{m+\beta}[ h]_{B(y,1)\cap A}   +     \|h \|_{m, A}
	\end{align*}
	Let us also define, for a function of class $C^{2,\alpha}(A)$,
	\be\label{star}
	\| \phi\|_{*,m-2, A} =  \|D^2\phi\|_{m,\beta, A}   +  \|D\phi\|_{m-1,A} + \|\phi\|_{m-2,A} .
	\ee
	In this notation we omit the dependence on $A$ when $A= \R^2$. The following is the main result of this section.

	\begin{prop}\label{prop1} (Proposition 6.1 in \cite{ddmw2}) There is $C>0$ such that for all sufficiently large $R$ and a differential operator $B_i$ as in $\equ{defB}$ with estimates \eqref{biest} and \eqref{By},
		Problem $\equ{eee}$ has a solution $\phi = T_i[h]$ for certain scalars $c_{ij}= c_{ij}[h]$, that defines a linear operator of $h$ and satisfies
		\begin{align*} 
			\| \phi\|_{*,m-2, B_R} \ \le\ C \|h \|_{m,\beta, B_R} .
		\end{align*}
		In addition, the linear functionals $c_i$ can be estimated as
		\begin{align*} 
			c_{i0}[h]\, = & \, \gamma_0\int_{B_{R}}  h Z_0  + O( R^{2-m})  \|h \|_{m,\beta, B_R} , \\    c_{ij}[h]\, = &\, \gamma_j\int_{B_{R}}  h Z_j  + O( \ve \mu_i \log|\log \ve|)  \|h \|_{m,\beta, B_R}, \ j=1,2.
		\end{align*}
		where
		$\gamma_j^{-1} = \int_{\R^2} e^{\Gamma_0} Z_j^2$, $j=0,1,2$.
	\end{prop}


	\section{Solving the inner-outer gluing system}\label{sette}
	
	
	We let $X^o$ be the Banach space of all functions $\psi \in C^{2,\beta}(\R^2)$ such that
	$$
	\|\psi \|_\infty < +\infty,
	$$
	and formulate the outer equation \equ{outer} as the fixed point problem in $X^o$,
	\begin{align*}
		\psi  =  \TT^o [   G(\psi, \phi)   ], \quad \psi\in X^o
	\end{align*}
	where ${\mathcal T}^o$ is defined in Proposition \ref{prop2}, while  $G$ is the operator given by \equ{GG}. 
	
	We formulate the projected problem \eqref{innerp} as the one of finding $(\phi_i , c_{ij}) $ where
	$$
	\phi_i = \phi_{i,1} + \phi_{2,i}
	$$
	with 
	\be\label{inner1}
	\Delta_y \phi_{i,1}  +    e^{\Gamma_0} \phi_{i,1} +   B_i[\phi_{i,1}]  + B_i[\phi_{i,2}]+ H_i(\phi, \psi)  =\sum_{j=0}^2 c_{ij} e^{\Gamma_0} Z_j   \inn B_R
	\ee
	where $Z_j$ are given by \eqref{defZj}, and
	\be\label{inner2}
	\Delta_y \phi_{i,2}  +    e^{\Gamma_0} \phi_{i,2} +  c_{i0} e^{\Gamma_0} Z_0  = 0 \inn \R^2 .
	\ee
	Problem \eqref{inner1} is formulated using the operator $T_i$ in Proposition \ref{prop1}, with 
	\begin{align*}
		c_{ij} &= c_{ij} [ H_i(\phi, \psi) +B_i(\phi_{i,2}) ], \quad j=0,1,2\\
		\phi_{i,1} &\in X_*, \quad   \phi_{i,1} =  T_i(  H_i(\phi, \psi) + B_i[\phi_{i,2}] ),
	\end{align*}
	where $X_*$ be the Banach space of functions $\phi \in C^{2, \beta} (B_R)$ such that
	$$
	\| \phi \|_{*,m-2,B_R} <\infty
	$$
	(see \eqref{star}). 
	
	\medskip
	Problem \eqref{inner2} is formulated using the operator $\mathcal{T}$ in Lemma \ref{lemat}
	$$
	\phi_{i,2} =  \TT [  c_{i0}[ H_i(\phi, \psi) +B_i(\phi_{i,2}) ]e^{\Gamma_0} Z_0   ], 
	$$
	in fact $\phi_{i,2}$ is a radial function satisfying
	\be\label{phi2ex}
	\phi_{i,2}(y) = c_{i0}[ H_i(\phi, \psi) +B_i(\phi_{i,2}) ]\left( \frac{4}{3}\frac{|y|^2-1}
	{|y|^2+1}\log(1+|y|^2)-\frac 83\frac{1}
	{|y|^2+1}\right).
	\ee
	
	Having in mind the a-priori bound in \eqref{cota}, \eqref{cotaa} in Lemma \ref{lemat}, it is natural  to ask that $\phi_{i,2} \in C^{2,\beta} (\R^2) $,
	\be\label{norm**}
	\begin{aligned}
		\| \phi \|_{**,\beta}&=\sup\limits_{y\in \RR^2}
		\frac{1}{\log (1+ |y|)} \Biggl[(1+|y|^{2+\beta})  [D^2_y \phi]_{B_1(y),\beta}  \\
		&+(1+|y|^2)  |D^2_y \phi (y)| +(1+|y|) | \nn \phi (y)|  +  | \phi (y)|\Biggl]
	\end{aligned}
	\ee
	and denote by $X_{**}$ the Banach space of functions $\phi \in C^{2,\beta}$ with $\| \phi \|_{**, \beta} <\infty$.
	
	\medskip
	\begin{prop}\label{ff}
		Let $\beta \in (0,1)$.
		There exist  positive constants $C$ and $\sigma_1>0$, functions $\psi \in X^o$, $\bar \phi_1 = (\phi_{1,1}, \ldots , \phi_{1,N}) \in X_*^N$, $\bar \phi_2 = (\phi_{2,1}, \ldots , \phi_{2,N}) \in X_{**}^N$, and constants $c_{ij}$,  $i=1, \ldots , N$, $j=1,2$, solutions to \eqref{outer}-\eqref{inner1}-\eqref{inner2} such that
		\be \label{estim}
		\begin{aligned}
			\| \psi \|_\infty &\leq C \ve^{1+\sigma_1} , \quad \| \phi_{i,1}\|_{*, m-2, B_R} \leq CR^{-1}, \\
			\, &\| \phi_{i,2}\|_{**, \beta} \leq C R^{1-m}  ,\: i:=1,\ldots,N\,
		\end{aligned}
		\ee
		with $R=2\delta_1/(\ve\mu_i|\log\ve|).$
	\end{prop}
	
	\medskip
	\proof
	Problem \eqref{innerp}-\eqref{outer} consists in 
	finding  $\psi$, $\bar \phi_{1}$, $\bar \phi_{2}$  solution of the fixed point problem
	\be \label{fix}
	(\psi, \bar \phi_{1} , \bar \phi_{2}) = {\mathcal A} (\psi, \bar \phi_{1} , \bar \phi_{2})
	\ee
	given by
	\begin{equation}\label{fixp}
		\begin{aligned}
			\psi  &=  \TT^o [   G(\psi, \bar \phi_{1} + \bar \phi_{2})   ], \quad \psi\in X^o\\
			\phi_{i,1} &=  T_i[  H_i(\bar \phi_{1} + \bar \phi_{2}, \psi) +B_i(\phi_{i,2})],\quad \phi_{i,1}\in X_*\\
			\phi_{i,2} &=  \TT [  c_{i0}[ H_i(\bar \phi_{1} + \bar \phi_{2}, \psi) +B_i(\phi_{i,2}) ]e^{\Gamma_0} Z_0   ],\quad \phi_{i,2}\in X_{**}.
		\end{aligned}
	\end{equation}
	
	Let  $m>2$ and define
	\be \label{ball}
	\begin{aligned}
		B_M&=\{ (\psi, \bar \phi_1 , \bar \phi_2  ) \in X^o\times X_*^N\times X_{**}^N  \, : \,  \\
		& \| \psi \|_\infty \leq M \ve^{1+\sigma_1}
		, \, \| \phi_{i,1}\|_{*, m-2, B_R} \leq M R^{-1} , \\
		&\, \| \phi_{i,2}\|_{**, \beta} \leq M  R^{1-m}
		,\: i:=1,\ldots,N\,  \},
	\end{aligned}
	\ee
	for some positive constant $M$ independent of $\ve$ and $\sigma_1$.
	We shall solve \eqref{fix}-\eqref{fixp} in $B_M$.
	
	We first show that ${\mathcal A} (B_M) \subset B_M$. Assume that  $(\psi, \bar \phi_1 , \bar \phi_2  ) \in B_M$. We first want to show that ${\mathcal A} (\psi, \bar \phi_1 , \bar \phi_2  ) \in B$.
	By definition of $X_*$ and $X_{**},$ we have 
	\begin{align*}
		A_i(\phi_{i})&\leq \frac C{1+|x|^\nu}\left(|\log \ve|^2|\phi_{i,1}+\phi_{i,2}|+\frac{|\log\ve|}{\ve \mu_i}|D_y (\phi_{i,1}+\phi_{i,2})|\right) \\ 
		&\leq \frac{R^{2-m}|\log\ve|^2}{1+|x|^\nu}\|\phi_{i,1}\|_{*,m-2,B_R}+\frac{|\log\ve|^3}{1+|x|^\nu} \|\phi_{i,2}\|_{**,\beta} 
	\end{align*}
	for some $\nu>2$ and for $x$ in a subset of $B((x_0,0),1),$ and $A_i(\phi_i)=0$ elsewhere.
	From \eqref{GG} and \eqref{house1}, we get that
	\be \label{f1g}
	\begin{aligned}
		|G(\psi , \bar \phi_{1} + \bar \phi_{2} )| &\leq {C \over 1+ |x|^\nu}  \ve^{1+b} \left(1+ \sum\limits_{i=1}^N|\phi_{i,1}+\phi_{i,2}|^2 \eta_i + |\psi|^2+ |\psi| \right) \\
		&+ {C \over 1+ |x|^\nu}
		|\log \ve|^2 \sum\limits_{i=1}^N(R^{2-m}\|\phi_{i,1}\|_{*, m-2, B_R} + |\log\ve|\| \phi_{i,2} \|_{**, \beta} ).
	\end{aligned}
	\ee
	From Proposition \ref{prop2}, we get that
	\be \label{f1}
	\begin{aligned}
		\| \psi \|_\infty &= \| {\mathcal T}^o ( G(\psi , \bar \phi_{1} + \bar \phi_{2}  ) ) \|_\infty \leq C \ve^{1+\sigma_1}
	\end{aligned}
	\ee
	where $\sigma_1=\min\{b,m-2-\sigma\}$ for $\sigma>0$ small.
	From \eqref{house2}, \eqref{ba1}, we get, for some $a\in (0,1)$,
	\be\label{f1b}
	\begin{aligned}
		|H_i(\bar\phi_1 + \bar\phi_2 , \psi ) | &\leq  |\tilde E_i|+
		{8 \over 1+ |y|^{2}} \left({1\over (1+ |y|^{2}) } +  {C \ve\mu_i \log|\log \ve | \over  (1+ |y|^{a})}\right)|\psi| \\
		&+ {C \over (1+ |y|)^4} \left(|\psi|^2 + \sum\limits_{i=1}^N(|\eta_i \phi_{i,1}|^2 + |\eta_i \phi_{i,2}|^2) \right).
	\end{aligned}
	\ee
	where
	$$
	|\tilde E_i| \leq C  { \ve\mu_i  \log|\log \ve | \over (1+ |y|^{2+a})} . 
	$$
	Using the assumptions on $\psi$, $\phi_{i,1}$ and $\phi_{i,2}$, we get that
	\be \label{f2}
	\|H_i(\bar \phi_1 + \bar \phi_2 , \psi ) \|_{m, \beta, B_R} \leq C \ve\mu_i \log|\log \ve| \leq CR^{-1}.
	\ee
	for $m<2+a$. From \eqref{defB}, \eqref{biest}, and \eqref{By}, we get
	\be\label{f2b}
	\begin{aligned}
		| B_i[\phi_{i,2}] | &\leq C\ve \mu_i \left(|D\phi_{i,2}|+|y| |D^2\phi_{i,2}| + {\log|\log \ve | \over 1+ |y|^{2+a}} |\phi_{i,2}|\right)
	\end{aligned}
	\ee
	and using again that $m<2+a$, we have
	\be \label{f3}
	\begin{aligned}
		\|B_i&(\phi_{i,2})  \|_{m, \beta, B_R} \leq C R^{m-2} \| \phi_{i,2} \|_{**, \beta}\\
		& + C R^{-1} \log|\log \ve |  \| \phi_{i,2} \|_{**, \beta} 
		\leq   C R^{m-2} \| \phi_{i,2} \|_{**, \beta} \leq CR^{-1}
	\end{aligned}
	\ee
	from Proposition \ref{prop1}, using estimates \eqref{f2} and \eqref{f3},  we find  that
	\be \label{f4}
	\| \phi_{i,1} \|_{*, m-2, B_R} \leq C R^{-1}
	\ee
	Now using \eqref{ee2}, \eqref{biest0} and the fact that $\phi_{i,2}$ is radial, we have 
	\begin{align*}
		\int\limits_{B_R}\tilde E_iZ_0=O(\ve^2\mu_i^2|\log\ve|), \quad \int\limits_{B_R}B_i(\phi_{i,2})Z_0=O( \| \phi_{i,2} \|_{**, \beta} R^2(\ve\mu_i)^2|\log\ve|).
	\end{align*}
	In addition,  since $\int_{B_R}e^{\Gamma_0(y)}Z_0=O(R^{-2})$ and  by regularity $\psi(x)=\psi(P_i)+\ve\mu_iA_iy\|\psi\|_\infty$ for $x\in B(P_i,\delta/|\log\ve|)$, we get
	\begin{align*}
		\int\limits_{B_R}[e^{\Gamma_0(y)}+b_i(y)]\psi Z_0=O(\ve\mu_i\|\psi\|_\infty),
	\end{align*}
	also we have $\int\limits_{B_R}{\mathcal N}_i(\sum\limits_{i=1}^N\eta_i \phi_i + \psi )Z_0=O( \sum\limits_{i=1}^N\|\phi_{i,1}\|^2_{*, m-2, B_R}),$ and consequently from the definition of $H_i$,
	$$
	\int\limits_{B_R} [ H_i(\phi, \psi) +B_i(\phi_{i,2}) ]Z_0= O\left(\frac{ \| \phi_{i,2} \|_{**, \beta}}{|\log\ve|}\right).
	$$
	By Proposition \ref{prop1}, we have 
	$$
	|c_{i0}[ H_i(\phi, \psi) +B_i(\phi_{i,2}) ]|\leq  
	C R^{2-m}\| H_i(\phi, \psi) +B_i(\phi_{i,2}) \|_{m,\beta,R} +O\left(\frac{\| \phi_{i,2} \|_{**, \beta}}{|\log\ve|}\right) 
	$$
	then 
	\be\label{f33}
	|c_{i0}[ H_i(\phi, \psi) +B_i(\phi_{i,2}) ]|\leq C  \| \phi_{i,2} \|_{**, \beta}
	\ee
	while from Lemma \ref{lemat} and \eqref{f33} we get
	\be \label{f5}
	\| \phi_{i,2} \|_{**, \beta} \leq C|c_{i0}[ H_i(\phi, \psi) +B_i(\phi_{i,2}) ]|\leq C R^{1-m}
	\ee
	Combining \eqref{f1}-\eqref{f4}-\eqref{f5}, we conclude that ${\mathcal A} (\psi, \phi_1 , \phi_2  ) \in B_M$ if we choose $M$ large enough (but independently of $\ve$) in the definition of the set $B_M$ in \eqref{ball}.
	
	We next show that ${\mathcal A}$ is a contraction map in $B_M$.
	Let $\varphi^j = \sum_{i=1}^N\eta_i (\phi_{i,1}^j + \phi_{i,2}^j) +\psi^j $,  for $j:=1,2$, such that
	$(\psi^j, \bar \phi_{1}^j, \bar \phi_{2}^j ) \in B_M$. Let $G(\varphi^j ) = G(\psi^j, \bar\phi_{1}^j + \bar \phi_{2}^j)$
	and observe that
	\be \label{estGG}\begin{aligned}
		\left| G(\varphi^1 ) - G(\varphi^2)  \right|&\leq |V (x) (\psi^1 - \psi^2)| \\
		& + \left(1-\sum\limits_{i=1}^N\eta_i \right) |N_{0 } (\varphi^1 ) - N_{0} (\varphi^2) |\\
		&+\sum\limits_{i=1}^N|A_i[\phi^1_{i,1}-\phi_{i,1}^2]| + \sum\limits_{i=1}^N|A_i[\phi_{i,2}^1 - \phi_{i,2}^2]|
	\end{aligned}
	\ee
	where the  terms are defined in \eqref{GG}.
	A direct computation gives that
	
	\begin{align*}
		&|V_{} (x) (\psi^1 - \psi^2)|+ (1-\sum\limits_{i=1}^N\eta_i ) | N_{0} (\varphi^1 ) - N_{0} (\varphi^2 ) | \\
		&\leq 
		C(\ve\mu)^2|\log \ve|^4
		\left( |\psi^1 - \psi^2| + \sum\limits_{i=1}^N\eta_i^2 (|\phi_{i,1}^1 - \phi_{i,1}^2|^2 + |\phi_{i,2}^1 - \phi_{i,2}^2|^2 ) \right)
	\end{align*}
	and as in \eqref{f1g}, 
	\begin{align*}
		|A_i[\phi_{i,1}^1-\phi_{i,1}^2]|& \leq  C \left(|L (\eta_i) |\phi_{i,1}^1 - \phi_{i,1}^2| + |K_{\ell j} \pp_{x_\ell} \eta_i \pp_{x_j} (\phi_{i,1}^1 - \phi_{i,1}^2) | \right) \\
		&\leq { C \ve^{m-2-\sigma} \over 1+ |x|^\nu } \| \phi_{i,1}^1 - \phi_{i,1}^2 \|_{*, m-2, B_R} .
	\end{align*}
	for $\sigma>0$ small.
	In order to estimate $A_i[\phi_2^1 - \phi_2^2]$, we observe that
	\begin{align*}
		\Delta_y & [\phi_{i,2}^1 - \phi_{i,2}^2]  +    f'(\Gamma_0) [\phi_{i,2}^1 - \phi_{i,2}^2] + c_{i0}^{12}   e^{\Gamma_0} Z_0  = 0 \inn \R^2  
	\end{align*}
	where
	\begin{align*}
		c_{i0}^{12}=   c_{i0} [H_i(\bar \phi_{1}^1+\bar\phi_{2}^1, \psi^1) +B_i(\phi_{i,2}^1) ] - c_{i0} [H_i(\bar \phi_{1}^2+\bar \phi_{2}^2, \psi^2) +B_i(\phi_{i,2}^2 )] .
	\end{align*}
	By definition,
	\begin{align*}
		&c_{i0}^{12}= \int_{B_R} \left[ B_i[\phi_{i,2}^1 - \phi_{i,2}^2]  + (e^{\Gamma_0(y)}+b_0(y))(\psi^1-\psi^2)  \right. \\
		&\mathcal N_i \left(\psi^1 + \sum_i\eta_i (\phi_{i,1}^1+ \phi_{i,2}^1)  \right) 
		\left. - \mathcal N_i \left(\psi^2 + \sum_i\eta_i (\phi_{i,1}^2+ \phi_{i,2}^2)  \right) \right] Z_0\,dy.
	\end{align*}
	Using \eqref{defB} and \eqref{house2}, we get
	\begin{align*}
		&|c_{i0}^{12} | \leq C\Biggl[{\frac{1}{|\log\ve|} } 
		\| \phi_{i,2}^1 - \phi_{i,2}^2\|_{**, \beta} +  \| \phi_{i,2}^1 - \phi_{i,2}^2\|_{**, \beta}^2
		+   \| \phi_{i,1}^1 - \phi_{i,1}^2\|_{*,m-2, \beta}^2 \\
		&+   \| \psi^1 - \psi^2\|_{\infty}^2 +R^{-1}[\| \psi^1 - \psi^2\|_{\infty}+ \| \phi_{i,1}^1 - \phi_{i,1}^2\|_{*,m-2, \beta} + \| \phi_{i,2}^1 - \phi_{i,2}^2\|_{**, \beta}]\Biggl].
	\end{align*}
	Then by \eqref{phi2ex} and \eqref{norm**}, we obtain
	\be\label{phi2cont}
	\| \phi_{i,2}^1 - \phi_{i,2}^2\|_{**, \beta} 
	\leq C |c_{i0}^{12} |.
	\ee
	Now we have 
	\begin{align*}
		|A_i[\phi_{i,2}^1-\phi_{i,2}^2]|& \leq  C \left(|L (\eta_i ) |\phi_{i,2}^1 - \phi_{i,2}^2| + |K_{\ell j} \pp_{x_\ell} \eta_i \pp_{x_j} (\phi_{i,2}^1 - \phi_{i,2}^2) | \right) \\
		&\leq { C |\log \ve|^3 \over 1+ |x|^\nu } \| \phi_{i,2}^1 - \phi_{i,2}^2 \|_{**, \beta} .
	\end{align*}
	Combining all these estimates in \eqref{estGG}, we obtain iterating \eqref{phi2cont}, that
	\begin{align*}
		\left| G(\varphi^1) - G(\varphi^2 )  \right|&\leq
		{C\ve^\sigma_0 \over 1+|x|^\nu} \Bigl( \| \psi^1-\psi^2\|_\infty + \sum\limits_{i=1}^N\| \phi_{i,1}^1 - \phi_{i,1}^2 \|_{*,m-2, B_R} \Bigl) \\
		& +{C \over 1+|x|^\nu} \frac{1}{|\log\ve|}\sum\limits_{i=1}^N \| \phi_{i,2}^1 - \phi_{i,2}^2 \|_{**,\beta}. 
	\end{align*}
	From Proposition \ref{prop2}, we have
	\be\label{psicont}
	\begin{aligned}
		\| \TT^o [G(\varphi^1 )] - \TT^o [G(\varphi^2 )]  \|_\infty &\leq
		C\ve^{\sigma_0}  \Bigl( \| \psi^1-\psi^2\|_\infty + \sum\limits_{i=1}^N \| \phi_{i,1}^1 - \phi_{i,1}^2 \|_{*,m-2, B_R} \Bigl)\\
		& +  \frac{C}{|\log\ve|}\sum\limits_{i=1}^N \| \phi_{i,2}^1 - \phi_{i,2}^2 \|_{**,\beta} ,
	\end{aligned}
	\ee
	for some $\sigma_0 >0$. Now let $T_i(\varphi^j)=T_i[  H_i(\bar \phi_{1}^j + \bar \phi_{2}^j, \psi^j) +B_i(\phi_{i,2}^j)]$  for $j=1,2$,
	then by Proposition \ref{prop1},
	\begin{align*}
		\|T_i&(\varphi^1)-T_i(\varphi^2)\|_{*,m-2,B_R}\leq C\Bigl[\|H_i(\bar\phi_1^1+\bar\phi_2^1,\psi^1)-H_i(\bar\phi_1^2+\bar\phi_2^2,\psi^2)\|_{m,\beta,B_R}\\
		&\qquad \qquad \qquad \qquad+\|B_i(\phi_{i,2}^1-\phi_{i,2}^2)\|_{m,\beta,B_R}\Bigl]\\
		&\leq  {C } \Bigl[
		\| \psi^1 - \psi^2\|_{\infty} +   \| \phi_{i,2}^1 - \phi_{i,2}^2\|_{**, \beta}^2
		+   \| \phi_{i,1}^1 - \phi_{i,1}^2\|_{*,m-2, \beta}^2+  \| \psi^1 - \psi^2\|_{\infty}^2 \\
		&+R^{-1}[\| \psi^1 - \psi^2\|_{\infty}+ \| \phi_{i,1}^1 - \phi_{i,1}^2\|_{*,m-2, \beta} + \| \phi_{i,2}^1 - \phi_{i,2}^2\|_{**, \beta}]\\
		& + R^{m-2} 
		\| \phi_{i,2}^1 - \phi_{i,2}^2\|_{**, \beta}\Bigl].
	\end{align*}
	As a consequence using \eqref{phi2cont} and \eqref{psicont},  we get that ${\mathcal A}$ is a contraction mapping in $B_M$ and  Problem \eqref{fix}-\eqref{fixp} has a fixed point. 
	\qed
	
	\section{The reduced problem}\label{otto}
	
	In Section \ref{sette} we proved the existence of a solution $(\phi_1, \ldots , \phi_N, \psi)$  to the coupled system of equations
	$$
	\Delta_y \phi_i +f'(\Gamma_0)\phi_i + B_i(\phi_i)+H_i(\phi,\psi)=\sum_{j=1}^2c_{ij}e^{\Gamma_0(y)}Z_j\quad \mbox{in}\quad B_R,
	$$
	for $i=1\ldots N$, and 
	$$
	L\psi+G(\psi,\phi)=0\quad \mbox{in}\quad \RR^2.
	$$
	The solution is described in Proposition \eqref{ff}, and estimates are contained in \eqref{estim}.

	In order to obtain an actual solution to our main Problem \eqref{ecuacion}, we need to show that the reduced system
	$$
	c_{ij}=c_{ij}[B_i(\phi_{i,2})+H_i(\phi,\psi)]=0\quad \mbox{for}\quad i=1\ldots N,\quad j=1,2
	$$
	can be solved provided the points $P_1, \ldots , P_N$ in \eqref{points}-\eqref{points1} are chosen properly.
	From Section \ref{sette} we get that
	$\| B_i(\phi_{i,2})+H_i(\phi,\psi)\|_{m, \beta , B_R} \lesssim \ve \mu \log |\log \ve|$. Hence 
	from Proposition \ref{prop1}, we obtain that
	$$
	c_{ij}=\gamma_j\int\limits_{B_R}[H_i(\phi,\psi)+ B_i(\phi_{i,2})]Z_j\, dy+\ve^{1+\sigma} {\bf Y} (P) ,
	$$
	for some $\sigma >0$. Here, and in the rest of this section, with ${\bf Y} (P)$ we denote a smooth function,  uniformly bounded  as $\ve \to 0$  for points  $P= (P_1, \ldots , P_N) $   satisfying  \eqref{points}-\eqref{points1}. The specific expression of this function changes from line to line, and even in the same line.
	
	Besides, since $H_i(\phi,\psi)= {\mathcal N}_i (\sum\limits_{i=1}^N\eta_i \phi_i + \psi ) +   \ttt  E_i +  (e^{\Gamma_0} + b_i ) \psi$, using the estimates \eqref{f1b} and \eqref{f2b}, we find
	$$
	\int\limits_{B_R}[{\mathcal N}_i (\sum\limits_{i=1}^N\eta_i \phi_i + \psi ) +     (e^{\Gamma_0} + b_i ) \psi+ B_i(\phi_{i,2})]Z_j
	=\ve^{1+\sigma} {\bf Y} (P).
	$$
	This fact implies that solving the reduced system 
	$c_{ij}=0$ is equivalent to prove  
	\be\label{condE}
	\int\limits_{B_R}\ttt E_iZ_j\,dy = \ve^{1+\sigma} {\bf Y} (P) \quad \mbox{for}\quad i=1,\ldots,N\quad \mbox{and}\quad j=1,2.
	\ee
	Formula \eqref{ee2} gives a rather explicit expression for  $\ttt E_i$, which is used to get
	$$
	\int\limits_{B_R}\ttt E_iZ_1\,dy=\ve\mu_i \kappa_i [M F_{1,i}(P)+\log(|\log\ve|)M{\bf Y}(P) 
	+  G_{1,i}(P)]
	$$
	where $M=\int\limits_{\RR^2}U(y)y_1Z_1(y)\,dy$,
	$$
	F_{1,i}(P)=-\left[\log \ve  \,   \left(  2  {R_i \over  \sqrt{(1+ R_i^2)^3} }   \, 
	-\alpha  {R_i \over \kappa_i\sqrt{1+ R_i^2} } \right) \right.
	\left.+ \sum_{j\not= i} {\kappa_j \over \kappa_i} 4 {[A_j^{-1} (P_i - P_j) ]_1  \over |A_j^{-1} (P_i - P_j)|^2} \right] 
	$$
	and
	$$
	G_{1,i}= c_{1,i}\int\limits_{\RR^2}U\Gamma_0y_1Z_1\,dy +
	{ R_i (3h^2+ R_i^2) \over 2h(h^2+ R_i^2)^{3\over 2} } M+  \ve\mu_i {\bf Y} (P) .
	$$
	On the other hand
	$$
	\int\limits_{B_R}\ttt E_iZ_2\,dy=\ve \mu_i \kappa_i M[F_{2,i}(P)+ {\bf Y}(P)  +\ve \mu_i {\bf Y} (P)]
	$$
	where
	$$
	F_{2,i}=-\left[  \sum_{j\not= i} {\kappa_j \over \kappa_i} 4 {[A_j^{-1} (P_i - P_j) ]_2  \over |A_j^{-1} (P_i - P_j)|^2}\right] .
	$$
	We recall now the form of the points $P_1 , \ldots , P_N$ as in \eqref{points}-\eqref{points1}:
	$$
	P_i= (a_i , b_i ) =(r_0+s,0)+
	\frac{(\hat a_i,\hat b_i)}{|\log\ve|}, \quad R_i = \sqrt{a_i^2 + b_i^2}
	$$
	and define
	$$\bar P_i=(\sqrt{h^2+r_0^2}\,\hat a_i,\hat b_i).
	$$
	Inserting this information in \eqref{condE} we obtain that 
	the reduced problem is
	\begin{align}
		\label{reduced1}
		\,   \left(  2  {hR_i \over  \sqrt{(h^2+ R_i^2)^3} }   \, 
		-\alpha  {hR_i \over \kappa_i\sqrt{h^2+ R_i^2} } \right) 
		& + 4\sum_{j\not= i} {\kappa_j \over \kappa_i}  {[(\bar P_i - \bar P_j) ]_1  \over |(\bar P_i - \bar P_j)|^2}\\
		&=\frac{\log|\log\ve|}{|\log\ve|}{\bf Y}(P) \nonumber \\
		\label{reduced2}
		\sum_{j\not= i} {\kappa_j \over \kappa_i} 4 {[(\bar P_i - \bar P_j) ]_2  \over |(\bar P_i - \bar P_j)|^2}=
		&\frac{{\bf Y} (P) }{|\log\ve|}
	\end{align}
	where again ${\bf Y}(P)$  denote a generic smooth function,  uniformly bounded  as $\ve \to 0$  for points  $P= (P_1, \ldots , P_N) $   satisfying  \eqref{points}-\eqref{points1}. 
	
	\medskip
	The non-linear system \eqref{reduced1}-\eqref{reduced2} is a perturbation of the 
	following limit problem \eqref{sysabintro0}, which for convenience we write using the complex notation
	\be \label{sysab2}
	\sum_{j\not= i}   {\kappa_j  \over {\bf P}_i - {\bf P}_j} =\left(   \kappa_i {h r_0 \over 2\sqrt{(h^2+ r_0^2)^3} }   \, 
	-\alpha  {h r_0 \over 4\sqrt{h^2+ r_0^2} } \right) 
	\ee
	for $i:=1,\ldots,N. $ Here ${\bf P}_j=({\bf P}_{j,1},{\bf P}_{j,2})$  is identified with the complex number ${\bf P}_j={\bf P}_{j,1} + i \, {\bf P}_{j,2}$.
	
	For all $i=1, \ldots , N$, let  $\F_i\colon \C^N\mapsto \C$ be the $i$-th left-hand side in \eqref{sysab2}, that is, 
	$$\F_i({\bf P})=\sum_{j\not= i}   {\kappa_j  \over {\bf P}_i - {\bf P}_j}
	\quad\mbox{for}\quad i:=1,\ldots,N,
	$$
	and let $U_i$ denote the right-hand side of \eqref{sysab2}
	$$U_i(r_0)= \left(   \kappa_i {h r_0 \over 2\sqrt{(h^2+ r_0^2)^3} }   \, 
	-\alpha  {h r_0 \over 4\sqrt{h^2+ r_0^2} }\right) .$$
	The point ${\bf P}^0$ satisfies
	$$
	\F_i({\bf P}^0)= U_i(r_0).
	$$
	We can calculate explicitly the derivative of $\F$ at ${\bf P}^0$, and we get
	$$
	d\F_{{\bf P}^0}=\left(\begin{array}{cccc}
		-\sum_{i=2}^N\kappa_iT_{1i} & \kappa_2T_{12} &...& \kappa_NT_{1N} \\
		\kappa_1T_{21} &-\sum_{i=1,i\neq 2}^N\kappa_iT_{2i} & ...& \kappa_NT_{2N} \\
		.... \\
		\kappa_1T_{N1} & \kappa_2T_{N2} & ... &
		-\sum_{i=1}^{N-1}\kappa_iT_{Ni} 
	\end{array}\right)
	$$
	where $T_{ij}=1/({\bf P}^0_i-{\bf P}_j^0)^2=T_{ji}$. A direct inspection gives that the vector  $\bf e_0=(1,\ldots,1) \in \C^N$ is an element of the kernel of $d\F_{{\bf P}^0}$. The non-degeneracy assumption on the point ${\bf P}^0$ means precisely that this is the only element in the kernel.

	\medskip
	We look for a solution to \eqref{reduced1}-\eqref{reduced2} as a small perturbation of ${\bf P}^0$. 
	Let ${\bf q}=(q_1,...,q_N)\in \C^N$, and redefine $\bar P_j$ in complex variables, we write 
	\begin{align*}
		\bar P_j &= {\bf P}_j^0+q_j,\quad  j=1,\ldots,N, 
	\end{align*}
	then the reduced problem \eqref{reduced1}-\eqref{reduced2} can be written as \be\label{reduced0}\F_i(\bar P)=U_i(R_i)+\bar\sigma_i\quad i=1,\ldots,N\ee
	with ${\rm Re}(\bar\sigma_i)=\frac{\log|\log\ve|}{|\log\ve|}\, {\bf Y} (P)$ and ${\rm Im}(\bar\sigma_i)=O\frac{1}{|\log\ve|}\, {\bf Y} (P)$.
	We have the expansions
	$$
	\F(\bar P)=\F({\bf P}^0)+d\F_{\bf P^0 }({\bf q})+O(|\bf q|^2)
	$$
	and
	\begin{align*}
		U_i(R_i)&=U_i(r_0)+U_i'(r_0)\left(s+\frac{{\rm Re}({\bf P}_{i}^0+q_{i})}{\sqrt{r_0^2+h^2}|\log\ve|}\right)\\&+O\left((\frac{|\bf q|}{|\log\ve|}+|s|)^2\right)+O\left(\frac{|s|}{|\log\ve|}\right)
	\end{align*}
	Thus \eqref{reduced0} takes the form
	\be \label{redexp}
	d\F_{\bf P^0 }({\bf q})={\GG}(s,{\bf q})+\frac{sh}{2}\left[\frac{h^2-2r_0^2}{(h^2+r_0^2)^{5 \over 2}}{\bf e}_1-\frac{\alpha}{2}\frac{h^2}{(h^2+r_0^2)^{3 \over 2}}{\bf e}_0\right]
	\ee
	where
	$
	\bf e_1=(\kappa_1,\ldots,\kappa_N)
	$
	and ${\rm Re}\, {\GG}(s,{\bf q})=O\left(\frac{\log|\log\ve|}{|\log\ve|}\right)$ as $\ve \to 0$ with higher order dependence in $s$ and ${\bf q}$. Also ${\rm Im}\, {\GG}(s,{\bf q})=O\left(\frac{1}{|\log\ve|}\right)$ as $\ve \to 0$. 
	Since $d\F$ 
	has a one dimensional kernel, we have the kernel of $(d\F)^T$ is also one dimensional and we found that is spanned by $\bf e_1$. From the value of $\alpha$ given by \eqref{alphadef}, the projection of the right hand side of \eqref{redexp}, onto $\bf e_1$ is equal to 
	${\GG}\cdot {\bf e}_1 - sh \frac{r_0^2\sum \kappa_i^2}{(h^2+r_0^2)^{5\over 2}}  $. Now we consider the following projected problem
	$$
	\begin{aligned}
		d\F_{\bf P^0 }({\bf q})=&{\GG}(s,{\bf q})
		+\frac{sh}{2}\left[\frac{h^2-2r_0^2}{(h^2+r_0^2)^{5 \over 2}}{\bf e}_1-\frac{\sum \kappa_i^2}{\sum \kappa_i}\frac{h^2}{(h^2+r_0^2)^{5 \over 2}}{\bf e}_0\right] \\
		&- \frac{{\GG}\cdot {\bf e}_1 - sh \frac{r_0^2\sum \kappa_i^2}{(h^2+r_0^2)^{5\over 2}}}{\sum \kappa_i^2} {\bf e}_1 
	\end{aligned}
	$$
	For each $s$, since ${\bf P}^0$ is a non-degenerate solution, we can solve this problem  for $q:=q(s)$. From the estimates on ${\mathcal G}$, we get that $|q|\lesssim {\log |\log \ve | \over |\log \ve |}$. By a fixed point argument we can find a solution $s=s^*$ with ${\GG}\cdot {\bf e}_1 - sh \frac{r_0^2\sum \kappa_i^2}{(h^2+r_0^2)^{5\over 2}}=0$. Since ${\rm Re}\, {\GG}(s,{\bf q})=O\left(\frac{\log|\log\ve|}{|\log\ve|}\right)$ as $\ve \to 0$, one has that $|s^*| \lesssim {\log |\log \ve | \over |\log \ve |}$. We thus get a solution of \eqref{redexp} 
	$$
	P_i=\left(\begin{array}{l}r_0+s^*+\dfrac{{\rm Re}({\bf P}^0_i+q_i(s^*))}{|\log \ve|\sqrt{h^2+r_0^2}}\\\dfrac{{\rm Im}({\bf P}^0_i+q_i(s^*))}{|\log \ve|}
	\end{array}\right)
	$$
	with the expected estimates. 
	This concludes the proof of our result.

	\medskip
	\section{Appendix} \label{appe1}
	
	\proof[Proof of Proposition \ref{prop2}]
	To solve Equation \equ{louter}
	we decompose $g$ and $\psi$ into Fourier modes as
	$$ g(x) = \sum_{j=-\infty}^{\infty} g_j(r) e^{ji\theta}, \quad    \psi(x) = \sum_{j=-\infty}^{\infty}  \psi_j (r) e^{ji\theta},\quad x=re^{i\theta} .$$
	Then
	$$
	L[\psi]  =  {1\over h^2 + r^2} \left( {h^2 \over r^2} +1 \right) \pp_\theta^2 \psi + {h^2 \over r} \pp_r \left( {r \over h^2 + r^2} \pp_r\psi \right).
	$$
	Thus this operator decouples the Fourier modes: equation
	\equ{louter}  becomes  equivalent to the following infinite set of ODEs:
	\be\label{LLk} \begin{aligned}
		&L_k[\psi_k] +  g_k(r) = 0 , \quad r\in (0,\infty),
		\\
		&L_k[\psi_k]:=
		{h^2\over r} \big( {r \over r^2 + h^2} \psi_k'\big)'  -  {k^2\over h^2 + r^2} \left( {h^2 \over r^2} +1 \right)\psi_k   , \quad k\in \Z.
	\end{aligned}  \ee
	The operator $L_k $ when $r\to 0$ or $r\to +\infty$ resembles
	$$\begin{aligned}
		L_k[p]\ \sim &\ \frac 1r(rp')' - \frac {k^2p} {r^2} \ass r\to 0 \\
		L_k[p]\ \sim &\ \frac{ h^2}{r}(\frac 1r p')' - \frac {k^2p} {r^2} \ass r\to +\infty
	\end{aligned}
	$$
	For $k\ge 1$, $L_k$ satisfies the maximum principle. This gives the existence of a positive function $z_k(r)$ with $L_k[z_k]=0$ with
	$$
	z_k(r) \sim  r^k \ass r\to 0, \quad 
	z_k(r)\sim   r^{\frac 12 } e^{(k/h)r} \ass r\to +\infty.
	$$
	Take $k=1$.  The function $$
	\bar \psi(r) =   z_1(r) \int_r^\infty \frac {(1+ s^2) ds}{ s z_1(s)^2} \int_0^s \frac 1{1+ \rho^\nu } z_1(\rho) \rho\, d\rho  
	$$
	solves $ L_1[\bar \psi]  + \frac 1{1+ r^\nu } = 0$, and satisfies the bounds
	$$\begin{aligned}
		\bar \psi(r) =O(  r^2 ) \ass r\to 0 \\
		\bar \psi(r) = O( r^{-\nu+2}) \ass r\to +\infty.
	\end{aligned}
	$$
	We take $\bar \psi$ as a barrier for the equation at $k=1$. Besides, this functions works as a barrier also for
	$k\ge 2$.
	
	For $k\geq 2$, the function
	$$
	\psi_k(r) =  z_k(r) \int_r^\infty \frac {(1+ s^2 ) ds}{ s z_k(s)^2} \int_0^s h_k(\rho) z_k(\rho) \rho\, d\rho  ,
	$$
	is the unique decaying decaying solution   \equ{LLk},  and it satisfies the estimate
	$$
	|\psi_k(r)| \le  \frac 4{k^2}  {\|g\|_{\nu}} \bar \psi(r) .
	$$
	since
	$$
	|g_k(r)|  \le \frac {\|g\|_{\nu}} {1+ r^\nu  }.
	$$
	If $k=0$, the solution is given by the explicit formula
	$$
	\psi_0(r)  = -\int_0^r  \frac{1 + s^2}{h^2  s}
	\, ds \int_0^s  g_0(\rho) \rho \, d\rho  
	$$
	and satisfies the bound
	$$
	|\psi_0(r)| \le      C\|g\|_{\nu} (1+  r^2) .
	$$
	The function
	$$
	\psi(x) :=  \sum_{j=-\infty}^{\infty}  \psi_j (r) e^{ji\theta},
	$$
	with the $\psi_k$ being the functions built above, clearly  defines a linear operator of $g$  and satisfies estimate \equ{estimate}.
	The proof is concluded. \qed
	
	\medskip
	
	\proof[Proof of Lemma \ref{lemat}]
	Setting $y=re^{i\theta}$, 
	we write
	$$ h(y) =  \sum_{k=-\infty}^{\infty}  h_k (r) e^{ik\theta} , \quad \phi(y)= \sum_{k=-\infty}^{\infty} \phi_k (r) e^{ik\theta}
	$$
	The equation is equivalent to
	\be
	\LL_k [\phi_k] + h_k(r) = 0 , \quad r\in(0,\infty)
	\label{01}\ee
	where
	$$
	\LL_k [\phi_k] =   \phi_k''  + \frac 1r \phi_k'  +  e^{\Gamma_0}  \phi_k   - \frac{k^2}{r^2} \phi_k
	$$
	Using the formula of variation of parameters,  the following formula (continuously extended to $r= 1$ defines a smooth solution of \equ{01} for $k=0$:
	$$
	\phi_0  (r) =     -z(r)\int_{1}^r  \frac {ds}{ s z(s)^2}  \int_0^s h_0(\rho) z(\rho) \rho\,ds, \quad z(r) = \frac{r^2 - 1}{ 1 + r^2}
	$$
	Noting that $\int_0^\infty h_0(\rho) z(\rho) \rho\,ds = \frac 1{2\pi}\int_{\R^2} h(y)Z_0(y)\, dy      $
	we see that this function satisfies
	$$
	|\phi_0(r)| \, \le \,   C\big[ \, \log (2 + r) \Big| \int_{\R^2} h(y)Z_0(y)\, dy      \Big|   \, +\,  (1+r)^{2-m} \|h\|_{ m} \big ].
	$$
	Now we observe that
	$$
	\phi_k  (r) =     -z(r)\int_0^r  \frac {ds}{ s z(s)^2}  \int_0^s h_k(\rho) z(\rho) \rho\,ds, \quad z(r) = \frac{4r}{ 1 + r^2}
	$$
	solves \equ{01} for $k =-1,1$ and satisfies
	$$
	|\phi_k(r)| \, \le \,   C\big[ \,  (1+ r) \sum_{j=1}^2 \Big| \int_{\R^2} h(y)Z_j(y)\, dy      \Big|   \, +\,  (1+r)^{2-m} \|h\|_{ m} \big ].
	$$
	For $k=2$ there is a function $z(r)$ such that $\mathcal L_2[z] = 0 $,  $z(r)\sim r^{2}$ as $r\to 0$ and as $r\to \infty$.  For $|k|\ge 2$ we have that
	$$
	\bar\phi_k  (r) =     \frac 4{ k^2} z(r)\int_0^r  \frac {ds}{ s z(s)^2}  \int_0^s |h_k(\rho)| z(\rho) \rho\,ds,
	$$
	is a positive supersolution for equation \equ{01}, hence the equation has a unique solution $\phi_k$  with $|\phi_k(r)| \le \bar\phi_k  (r)$. Thus
	$$
	|\phi_k(r)| \, \le \,    \frac C {k^2}   (1+r)^{2-m} \|h\|_{ m}, \quad |k|\ge 2.
	$$
	Thus $\phi_k$ defined  $$\phi(y)= \sum_{k=-\infty}^{\infty} \phi_k (r) e^{ik\theta} $$  defines a linear operator of functions $h$ which is a solution of equation \equ{00} which, adding up the individual estimates above,  it satisfies the estimate \eqref{cota}.
	As a corollary we find that similar bounds are obtained for first and second derivatives. In fact, let us set
	for a large $y = R e$, $R=|y|\gg 1$,
	$
	\phi_R(z)  =  {R^{m-2}} \phi (R(e+z))    .
	$
	Then in a neighborhood of $y$, we find 
	$$
	\Delta_z\phi_R  +  \frac{8R^2}{(1+R^2|e+ z|^2)^2}\phi_R +  h_R (z)   = 0, \quad |z| <\frac 12
	$$
	where $h_R(z) = R^{m}h (R(e+z)) $.
	Let us set,
	$$
	\delta_i =  \Big| \int_{\R^2}  h Z_i \Big | , \quad i=0,1,2.
	$$
	Then  from \equ{cota}, and a standard elliptic estimate we find
	$$\|\nn_z \phi_R \|_{L^\infty( B_{\frac 14}(0) )} + \|\phi_R \|_{L^\infty( B_{\frac 12}(0) )} \, \le\,  C\Big [ \delta_0 R^{m-2}\log R +   \sum_{i=1}^2\delta_i R^{m-1} +   \| h\|_{m} \Big ],$$
	using that $ \|h_R \|_{L^\infty(B_{\frac 12} (0))} \ \le\  C\| h\|_{m}.$
	Now since    $[h_R ]_{B_{\frac 12} (0)} \ \le\  C\| h\|_{m,\alpha}$, from interior Schauder estimates and the bound for $\phi_R$  we then  find
	$$
	\|D^2_{z} \phi_R \|_{L^\infty( B_{\frac 14}(0) ) } +  [ D^2_{z} \phi_R ]_{B_{\frac 14}(0), \alpha} \, \le\,  C\Big [ \delta_0 R^{m-2}\log R +   \sum_{i=1}^2\delta_i R^{m-1} +   \| h\|_{m,\alpha} \Big ].
	$$
	From these relations, estimates \equ{cota} and \equ{cotaa} follow. \qed

	\medskip
	
	\proof[Proof of Proposition \ref{prop1}]
	We consider a standard linear extension operator $h\mapsto \ttt h $ to entire $\R^2$,
	in such a way that the support of $\ttt h$ is contained in $B_{2R}$ and $\|\ttt h\|_{m,\beta} \le C\|h\|_{m,\beta, B_R}$ with $C$ independent of all large $\bar R$. The operator $B_i $ as defined in \eqref{defB} and the coefficients are of class $C^1$ in entire $\R^2$ and  have compact support in $B_{2R}$. 
	Then we consider the auxiliary problem in entire space
	\be
	\Delta \phi  +  e^{\Gamma_0} \phi  + B_i[\phi]  + \ttt h(y)  = \sum_{j=0}^2  c_{ij} e^{\Gamma_0} Z_j \inn \R^2
	\label{001} \ee
	where, assuming that $\|h\|_m<+\infty$ and $\phi$ is of class $C^2$,
	$c_{ij} = c_{ij}[h,\phi]$ are the scalars defined as so that
	$$
	\gamma_i \int_{\R^2} (B_i[\phi]  + \ttt h(y))Z_j  =    c_{ij}  ,  \quad \gamma_j^{-1} = \int_{\R^2}  e^{\Gamma_0} Z_j^2
	$$
	For $j=1,2$, by \eqref{biest} and \eqref{By}, we have  $B_i[Z_j]=  O((1+|y|)^{-2})\ve\mu_i + O((1+|y|)^{-(3+a)})\ve\mu_i \log|\log \ve| $. Similarly for $j=0$,  we have  $B_i[Z_0]=  O((1+|y|)^{-3})\ve\mu_i + O((1+|y|)^{-(2+a)})\ve\mu_i \log|\log \ve| $. Since $m>2$, we get
	$$
	\int_{\R^2} B_i[\phi]Z_j =   \int_{\R^2} \phi \tilde B_i[Z_j]  = O(\|\phi\|_{m-2}) \ve \mu_i \log|\log \ve|   .
	$$
	where $\tilde B_i$ have same estimates as $B_i$ mentioned above. On the other hand
	$$
	\int_{\R^2\setminus B_R} h(y)Z_0 = O(R^{2-m})\|h\|_{m,\beta,B_R},\:\:
	\int_{\R^2\setminus B_R} h(y)Z_j = O(R^{1-m})\|h\|_{m,\beta,B_R}
	$$
	for $j=1,2.$ In addition, we readily check that
	$$ \| B_i[\phi] \|_{m,\beta}  \le  C \frac{\delta}{|\log\ve|}  \|\phi\|_{*, m-2,\beta},
	$$
	where
	$$
	\|\phi\|_{*, m-2,\beta} =  \|  D^2_y\phi \|_{m,\beta}   + \|  D_y\phi \|_{m-1}+ \|\phi \|_{m-2}.
	$$
	Let us consider the Banach space $X$ of all $C^{2,\beta}(\R^2)$
	functions with
	$
	\|\phi\|_{*, m-2,\beta} <+\infty.
	$
	We find a solution of \equ{001}
	if we solve the equation
	\be\label{fp}
	\phi  =   \mathcal A  [\phi]  +  \mathcal H ,\quad \phi \in X
	\ee
	where
	$$
	\mathcal A  [\phi]
	=  \mathcal T\Big [B_i[\phi] -\sum_{j=0}^2 c_{ij}[0, \phi] e^{\Gamma_0}Z_j  \Big] ,\quad
	\mathcal H  = \mathcal T\Big [  \ttt h   - \sum_{j=0}^2 c_{ij}[\ttt h,0] e^{\Gamma_0}Z_j  \Big] .
	$$
	and $\mathcal T$ is the operator built in Lemma  \ref{lemat}.
	We observe that
	$$
	\|\mathcal A  [\phi] \|_{*, m-2,\beta} \le C\frac{\delta}{|\log\ve|}  \|\phi\| _{*, m-2,\beta}, \quad \|\mathcal H \|_{*, m-2,\beta} \le C \|h \|_{ m,\beta, B_R}.
	$$
	So we find that Equation \equ{fp} has a unique solution, that defines a linear operator of $h$, and satisfies
	$$
	\|\phi \|_{*, m-2,\beta} \ \le\  C \|h \|_{ m,\beta, B_R}
	$$
	The result of the proposition follows by just setting $T_i[h] = \phi\big|_{B_R}$. The proof is concluded. \qed

	\medskip\noindent{\bf Acknowledgments:}
	Ignacio Guerra was supported by Proyecto Fondecyt Regular 1180628. 
	Monica Musso has been supported by EPSRC research Grant EP/T008458/1.

\end{document}